\documentclass[onecolumn,draftcls,11pt]{IEEEtran}
\usepackage{cite}
\usepackage{color}
\usepackage[pdftex]{graphicx}
\graphicspath{{fig/}{jpeg/}}
\usepackage[cmex10]{amsmath}
\usepackage{amssymb}
\usepackage{epstopdf}
\usepackage{multirow}
\usepackage{pifont}
\newcommand{\cmark}{\ding{51}}%
\newcommand{\xmark}{\ding{55}}%
\usepackage{algorithm}
\usepackage{tikz}
\usetikzlibrary{shapes,arrows,matrix}
\usetikzlibrary{decorations.pathmorphing} 
\usetikzlibrary{fit}					
\usetikzlibrary{backgrounds}	
\usepackage{verbatim}
\usepackage{algpseudocode}
\usepackage{amsmath,amssymb,lipsum}
\usepackage{hyperref}
\usepackage{comment}
\usepackage{graphicx}
\usepackage[font=small]{caption}
\usepackage{subcaption}
\usepackage{mathtools,enumerate}
\input{mysymbol.sty}
\usepackage{needspace}

\usepackage{theorem}
\newtheorem{thm}{Theorem}
\newtheorem{lemma}{Lemma}
\newtheorem{proposition}{Proposition}
\newtheorem{corollary}{Corollary}

\theoremstyle{definition}
\newtheorem{assumption}{Assumption}
\newtheorem{remark}{Remark}

\def \QED {$\blacksquare$}

\def \x {{\mathbf{x}}}

\def \y {{\mathbf{y}}}
\def \h {{\mathbf{h}}}

 \providecommand{\Ex}[1]{\mathbb{E}\left[#1\right]}
 \providecommand{\abs}[1]{\left|#1\right|}
 \providecommand{\norm}[1]{\left\|#1\right\|}
 \providecommand{\ip}[1]{\boldsymbol{\langle}#1\boldsymbol{\rangle}}

    \def \eps {\epsilon}
    \def \mbE {{\mathbb{E}}}
    \def \mbR {{\mathbb{R}}}
    \def \bbx {{\mathbf{x}}}
    \def \bby {{\mathbf{y}}}
    \def \ccalX {{\mathcal{X}}}
    
    \def \ccalL {{\mathcal{L}}}
    \def \cL {{\mathcal{L}}}
    \def \L {{\hat{\cL}}}
    \def \bbs {{\mathbf{s}}}
    
    \def \bbd {{\mathbf{d}}}
    \def \bbh {{\mathbf{h}}}
    \def \bbtheta {{\boldsymbol{\theta}}}
    \def \bblambda {{\boldsymbol{\lambda}}}
    \def \bblam {{\bblambda}}    
    \def \Rn {{\mathbb{R}}}
    \def \th {{\bbtheta}}
    \def \lam {{\bblambda}}
    
    \def \d {{\bbd}}
    \def \cX {{\ccalX}}
    \def \O {{\mathcal{O}}}
    \def \xu {\x_{\upsilon}}
    \def \H {{\mathbf{H}}}
    \def \px {{\mathcal{P}_{\mathcal{X}}}}
    
    \def \nx {{\nabla_\x}}
    \def \M {{\mathbf{M}}}
    \def \X {{\mathbf{X}}}
    
    \def \lus {{\lam_{\upsilon}^\star}}
    \def \xus {{\x_\upsilon^\star}}
    \def \F {{\mathcal{F}}}
    \def \sX {{\mathsf{X}}}
    \def \tx {{\tilde{\x}}}
    \def \one {{\textbf{1}}}
    \def \s {{\bbs}}
    \def \cI {{\mathcal{I}}}
    
  \tikzstyle{agent}=[circle,
  thick,
  minimum size=.5cm,
  draw=blue!80,
  fill=blue!20]

  \tikzstyle{neighbor}=[circle,
  thick,
  minimum size=.5cm,
  draw=cyan!85!black,
  fill=cyan!40,
  decorate,
  decoration={random steps,
  	segment length=2pt,
  	amplitude=2pt}]
  
  \tikzstyle{local_nat}=[rectangle,
  thick,
  minimum size=.5cm,
  draw=red!80,
  fill=red!20]
  
  \tikzstyle{glob_nat}=[rectangle,
  thick,
  minimum size=.5cm,
  draw=red!100,
  fill=red!40]      
  
  \tikzstyle{background}=[rectangle,
  fill=gray!10,
  inner sep=0.2cm,
  rounded corners=5mm]
  
  \tikzstyle{background2}=[rectangle,
  fill=green!20,
  inner sep=0.2cm,
  rounded corners=5mm]
{\tiny }
\title{\vspace{-0cm}{Conservative Stochastic Optimization with
		Expectation Constraints}}
\author{Zeeshan~Akhtar,$^\star$
	Amrit~Singh~Bedi,$^\dagger$
	and~Ketan~Rajawat$^\star$
	 \thanks{
	        Zeeshan Akhtar and K. Rajawat are with the Department of Electrical Engineering,
	        Indian Institute of Technology Kanpur, Kanpur 208016, India (e-mail: zeeshan@iitk.ac.in,
	        ketan@iitk.ac.in). A. S. Bedi is with the U.S. Army Research Laboratory, Adelphi, MD, USA. (e-mail: amrit0714@gmail.com). }}
\begin{document}
\maketitle

\begin{abstract}
    This paper considers stochastic convex optimization problems where the objective and constraint functions involve expectations with respect to the data indices or environmental variables, in addition to deterministic convex constraints on the domain of the variables. Although the setting is generic and arises in different machine learning applications, online and efficient approaches for solving such problems have not been widely studied. Since the underlying data distribution is unknown a priori, {a} closed-form solution is generally not available, and classical deterministic optimization paradigms are not applicable. Existing approaches towards solving these problems make use of stochastic gradients of the objective and constraints that arrive sequentially over iterations. State-of-the-art approaches, such as those using the saddle point framework, are able to ensure that the optimality gap as well as the constraint violation decay as  $\O\left(T^{-\frac{1}{2}}\right)$ where $T$ is the number of stochastic gradients. The domain constraints are assumed simple and handled via projection at every iteration. In this work, we propose a novel conservative stochastic optimization algorithm (CSOA) that achieves zero average constraint violation and  $\O\left(T^{-\frac{1}{2}}\right)$ optimality gap.
    
	Further, we also consider the scenario where carrying out a projection step onto the convex domain constraints at every iteration is not viable. Traditionally, the projection operation can be avoided by considering the conditional gradient or Frank-Wolfe (FW) variant of the algorithm. The state-of-the-art stochastic FW variants achieve an optimality gap of $\O\left(T^{-\frac{1}{3}}\right)$ after $T$ iterations, though these algorithms have not been applied to problems with functional expectation constraints. In this work, we propose the FW-CSOA algorithm that is not only projection-free but also achieves zero average constraint violation with $\O\left(T^{-\frac{1}{4}}\right)$ decay of the optimality gap. The efficacy of the proposed algorithms is tested on two relevant problems: fair classification and structured matrix completion. 
\end{abstract}

\vspace{-0mm}
\section{Introduction}\label{sec:intro}
We consider the following constrained stochastic convex optimization problem
\begin{align} \label{prob_form} 
	\x^\star =& \argmin_{\bbx \in \cX}  \mbE_{\bbtheta}[f(\bbx,\bbtheta)] \\
	& \text{\ s.t. }  \mbE_{\bbtheta}[\bbh({\bbx, \bbtheta})]           \leq 0, \nonumber
\end{align}
where $\mathcal{X}\subseteq\mbR^{m}$ is a convex and compact set, $\th \in \Theta \subset \Rn^d$ is a random vector, and the functions $f:\Rn^m \times \Rn^d \rightarrow \Rn$ and $\h:\Rn^m \times \Rn^d \rightarrow \Rn^N$ are proper, closed, and convex in $\x$. The distribution of $\th$ is unknown and we only have access to a stochastic gradient oracle that, for a given $\x \in \ccalX$, outputs $(\nabla f(\x,\th), \h(\x,\th), \nabla \h(\x,\th))$ for some randomly sampled $\th \in \Theta$. Here, we use $\nabla$ without any subscript to mean gradient with respect to $\x$ when applied to either $f$ or $\h$. 

\subsection{Background}
The relatively simpler problem of minimizing $\mbE[f(\bbx, \bbtheta)]$ over $\x \in \cX$ has been widely studied within the stochastic approximation rubric. In most applications, such as those arising in machine learning \cite{mu2017stochastic}, finance \cite{li2017designing},  robotics \cite{derenick2007convex}, signal processing  \cite{koppel2018parallel}, and communications \cite{8255626}, the constraint set $\cX$ is easy to work with, giving rise to two classes of stochastic gradient descent (SGD) algorithms: projected SGD and Franke-Wolfe SGD. In projected SGD and its variants, the iterate is projected back into $\cX$ after taking a step in the negative stochastic gradient direction  \cite{lan2012optimal,nemirovski2009robust,johnson2013accelerating,parikh2014proximal}. Such algorithms are efficient if the projection operation is computationally cheap, as in the case of the box or simplicial constraints. Alternatively, conditional or Frank-Wolfe SGD approaches provide a projection-free alternative \cite{frank1956algorithm,hazan2012projection,hazan2016variance, mokhtari2018stochastic,Zee_ACC,lan2016conditional}. Here, each iteration requires minimizing a linear objective over $\cX$, which may be computationally cheaper in some cases, such as when $\cX$ is an $\ell_1$ or nuclear-norm ball \cite{jaggi2013revisiting}. For smooth functions, after $T$ calls to the stochastic gradient oracle, the optimality gap of the projected SGD decays as $\O(T^{-\frac{1}{2}})$  \cite{nemirovski2009robust} while that of the Frank-Wolfe SGD decays as $\O(T^{-\frac{1}{3}})$ \cite{mokhtari2018stochastic}.

The expectation-constrained problem in \eqref{prob_form} has not received as much attention, though it finds applications in machine learning \cite{zafar2015fairness}, signal processing \cite{nemirovski2006convex,chapelle2006semi}, communications  \cite{singh2019asynchronous}, finance \cite{rockafellar2000optimization,wang2008sample}, and control theory \cite{vzliobaite2017measuring}. In general, projected SGD \cite{lan2012optimal,nemirovski2009robust,johnson2013accelerating,parikh2014proximal} and its Frank-Wolfe extensions \cite{frank1956algorithm,hazan2012projection,hazan2016variance, mokhtari2018stochastic,Zee_ACC,lan2016conditional} are designed for simple constraints in $\cX$ and  cannot directly handle expectation constraints since the functional form of $\mbE_{\bbtheta}[\bbh({\bbx, \bbtheta})]$ is not known a priori but is only revealed via calls to the first order oracle. 

\begin{table*}[t]\centering
	\begin{tabular}{|c|c|c|c|c|}
		\hline
		References                     &  \begin{tabular}[c]{@{}c@{}} Projection\\ Free\end{tabular}    &  \begin{tabular}[c]{@{}c@{}} Expectation\\ Constraint\end{tabular}      &  \begin{tabular}[c]{@{}c@{}} Optimality\\ Gap\end{tabular} & \begin{tabular}[c]{@{}c@{}} Constraint\\ Violation\end{tabular} \\ \hline
		\cite{bedi2019asynchronous}		& \xmark 	& \cmark		& $\O\left(T^{-{1}/{2}}\right)$ & $\O\left(T^{-{1}/{4}}\right)$ \\
		G-OCO \cite{yu2017online}        & \xmark   & \cmark            & $\O\left(T^{-{1}/{2}}\right)$ & $\O\left(T^{-{1}/{2}}\right)$ \\
		\cite{madavan2019subgradient}     & \xmark      & \cmark     &$\O\left(T^{-{1}/{2}}\right)$ & $\O\left(T^{-{1}/{2}}\right)$\\
		\cite{mahdavi2012trading}     & \xmark      & \xmark     &$\O\left(T^{-{1}/{4}}\right)$ & \textbf{zero} \\
		$\textbf{CSOA}$  (This Paper)                     & \xmark        & \cmark    & $\O\left(T^{-{1}/{2}}\right)$ & \textbf{zero}\\ \hline
		SFW \cite{hazan2012projection}  
		& \cmark          & \xmark    & $\O\left(T^{-{1}/{4}}\right)                                                   $ & -                                                                                                   \\
		SVRF \cite{hazan2016variance}                     & \cmark       & \xmark     & $\O\left(T^{-{1}/{3}}\right)                                                   $ & -                                                    \\
		SCGD \cite{mokhtari2018stochastic}                   & \cmark         & \xmark   & $\O\left(T^{-{1}/{3}}\right)                                                   $ & -                                                     \\ $\textbf{FW-CSOA} $ (This Paper)                    & \cmark        & \cmark    & $\O\left(T^{-{1}/{4}}\right)$ & \textbf{zero} \\\hline \end{tabular}
	\caption{Summary of related works. All results are normalized to require $T$ calls to the stochastic gradient oracle.\vspace{-4mm}}
	\label{table1}
\end{table*} 
\textbf{Projection-based algorithms} for solving \eqref{prob_form} have been recently studied in  \cite{yu2017online,lan2020conditional,basu2019optimal,mahdavi2012stochastic,bedi2019asynchronous,zhang2019stochastic,madavan2019subgradient}.  These algorithms can be classified into three categories: offline, online, and stochastic, depending on the approach followed. \emph{Offline} algorithms are relevant for settings where  $f$ and $\h$ have finite-sum structure, and the full data set is available prior to the start of the algorithm. The cooperative stochastic approximation (\textbf{CSA}) algorithm proposed in  \cite{lan2020conditional} considered the special case of a single expectation constraint, while its extended version proposed in \cite{basu2019optimal} studied the case of multiple expectation constraints. In both approaches, the optimality gap as well as the constraint violation decay as $\O(T^{-\frac{1}{2}})$ after $T$ calls to the oracle. 

\emph{Stochastic} algorithms apply to the setting considered here with stochastic gradient oracle returning gradients calculated for independent identically distributed (i.i.d.) samples of $\th$. A proximal method of multipliers approach is proposed in \cite{zhang2019stochastic} yielding convergence rates of $\O(T^{-\frac{1}{4}})$ for both the optimality gap and the constraint violation. A special case of \eqref{prob_form} for distributed and asynchronous implementation was considered in \cite{singh2019asynchronous}, where the optimality gap decayed as $\O(T^{-\frac{1}{2}})$ while the constraint violation decayed as $\O(T^{-\frac{1}{4}})$. A generalized version of problem in \ref{prob_form}, where the objective and the constraints are non-linear functions of expectations, was studied in \cite{thomdapu2019optimal, thomdapu2020stochastic}, but a slower convergence rate was achieved.  More recently, the saddle point algorithm in \cite{madavan2019subgradient} achieves $\O(T^{-\frac{1}{2}})$ rate for both objective as well as the constraint violation. Recently, \cite{boob2019stochastic} proposed a stochastic algorithm called constraint extrapolation (ConEx) for solving convex 	functional constrained problems, which utilize linear approximations of the constraint functions to define an extrapolation step. The work in \cite{boob2019stochastic} is closest in spirit to our work since it is also a primal-dual algorithm and achieves optimal convergence rates without the assumption of bounded primal-dual multipliers. However, different from \cite{boob2019stochastic}, the proposed algorithm achieves zero average constraint violation with the best possible convergence rate for the optimality gap.

Finally, the \emph{online} version of \eqref{prob_form} has also been studied, where the oracle may return arbitrary gradients that need not correspond to i.i.d. samples of $\th$. In such settings, the performance metrics of interest are regret and the cumulative constraint violation. The generalized online convex optimization (\textbf{G-OCO}) algorithm proposed in \cite{yu2017online} yielded a regret and cumulative constraint violation of $\O(T^{\frac{1}{2}})$ after $T$ calls from the oracle. While the online setting is more general than the stochastic one considered here, in the special case when the gradients come  from i.i.d. samples, \textbf{G-OCO} makes the optimality gap and the constraint violation decay as $\O(T^{-\frac{1}{2}})$. With functional but deterministic constraints, it is known from \cite{mahdavi2012trading,singh2019asynchronous} that the average constraint violation can be forced to zero, albeit at a penalty to the convergence rate of the optimality gap, which now becomes $\O(T^{-\frac{1}{4}})$. 

\textbf{Projection-free} algorithms for solving \eqref{prob_form} have not been proposed in the literature. Existing projection-free algorithms are only applicable to problems without expectation constraints.  For the instance, a projection-free version of \cite{boob2019stochastic} capable of handling  deterministic functional constraints was proposed in \cite{lan2020conditional}. For stochastic objective case, the work in \cite{hazan2012projection} utilized  mini-batches to obtain a convergence rate of $\O(T^{-\frac{1}{4}})$ given $T$ calls to the stochastic gradient oracle. Since then, the convergence rates have been improved to $\O(T^{-\frac{1}{3}})$ \cite{hazan2016variance, mokhtari2018stochastic}. 
\subsection{Contributions}
This work puts forth the conservative stochastic optimization algorithm (\textbf{CSOA}) that after $T$ stochastic gradient oracle calls make the optimality gap decay as $\O(T^{-\frac{1}{2}})$ while yielding zero average constraint violation. The result is remarkable since it demonstrates that the effect of the constraints can be completely nullified while ensuring that the optimality gap decays at the same rate as in projected SGD. In order to achieve this rate, we apply the stochastic saddle point algorithm to the following \emph{conservative} version of \eqref{prob_form}: 
\begin{align}
	\xu^\star&=\argmin_{\bbx} F(\x) := \mbE[f(\bbx,\bbtheta)]\nonumber\\
	& \text{s.t.} \;\;H_i(\x) +\upsilon\leq 0 \;\; \forall \;  i \in \{1,2,\cdots,N\}
	\label{new_prob}
\end{align} 
where, $H_i(\x):=\mbE[h_i(\bbx,\bbtheta)]$ and $\upsilon$ is an an algorithm parameter. Intuitively, $F(\xu^\star)$ is $\O(\upsilon)$ away from the optimal value $F(\x^\star)$. However, since we expect the optimality gap to decay as $\O(T^{-\frac{1}{2}})$, the approximation error due to the use of \eqref{new_prob} does not dominate, provided that $\upsilon = \O(T^{-\frac{1}{2}})$. At the same time, {the} presence of $\upsilon$ can cancel out the $\O(T^{-\frac{1}{2}})$ constraint violation that the saddle point algorithm would otherwise incur. The proposed approach prevents the updates from being too aggressive with regards to the constraint violation, while not sacrificing optimality. 

The idea of using conservative constraints is general and can also be applied to projection-free algorithms. We propose the Frank-Wolfe \textbf{CSOA} algorithm, which also solves \eqref{new_prob} while carefully selecting the value of $\upsilon$. Additionally, the analysis for the \textbf{FW-}\textbf{CSOA} algorithm is almost entirely novel since there are no other projection-free algorithms handling expectation constraints. After $T$ calls to the stochastic gradient oracle, the proposed algorithm achieves an optimality gap of $\O(T^{-\frac{1}{4}})$ and zero average constraint violation. Table \ref{table1} summarizes the convergence results for the proposed algorithms and compares them with the related works in the literature. In order to compare all the algorithms fairly, we have assumed that each algorithm makes $T$ calls to the stochastic gradient oracle and compared the resulting convergence rates for the optimality gap and the constraint violation. 
%
%
%
%
In summary, the contributions of this paper are as follows:
\begin{itemize}
	\item We propose \textbf{CSOA} to solve \eqref{prob_form} and show that it achieves a convergence rate of $\O(T^{-\frac{1}{2}})$ for  optimality gap without any constraint violation. The proposed algorithm utilizes augmented Lagrangian in the updates as in \cite{bedi2019asynchronous}, utilizes tools  from \cite{madavan2019subgradient} to improve the convergence rate, and draws upon the idea of forcing zero average constraint violation from  \cite{mahdavi2012trading}. 
	
	\item We propose the first projection-free algorithm \textbf{FW-}\textbf{CSOA} to solve \eqref{prob_form} that achieves a convergence rate of $\O(T^{-\frac{1}{4}})$ for the optimality gap without any constraint violation. The proposed algorithm is essentially a stochastic Frank-Wolfe version of the saddle point algorithm but again uses the idea of forcing  zero average constraint violation from  \cite{mahdavi2012trading}. The idea to make {the} algorithm projection-free is utilized as in \cite{mokhtari2018stochastic} and uses the momentum based gradient tracking as presented in \cite{cutkosky2019momentum}. 
	
	\item The proposed algorithms are rigorously tested on problems of fair classification and structured matrix completion, and their efficacy is demonstrated. 
	
\end{itemize}

\section{Algorithm Development}\label{Algo_Dev}
In this section, we propose \textbf{CSOA}  and \textbf{FW-CSOA} for solving \eqref{prob_form}. As stated earlier, the distribution of $\th$ is not known, and instead, only calls to the stochastic gradient oracle are allowed. Given $\x$, the oracle returns $(\nabla f(\x,\th_t), \h(\x,\th_t), \nabla \h(\x,\th_t))$ for an i.i.d. $\th_t$. As discussed in Sec.\ref{sec:intro}, we will concentrate on solving the following conservative version  of the problem defined in \eqref{prob_form}:
\begin{align}
	\xu^\star&=\argmin_{\x \in \mathcal{X}} F(\x) := \mbE[f(\x,\th)]\nonumber\\
	& \text{s.t.} \;\;\H(\x) +\upsilon\textbf{1}\leq 0
	\label{conservative_prob}
\end{align} 
where $\H(\x) := \mbE[\h(\x,\th)]$. Associating dual variables $\lam \in \Rn^N$ with the constraints, the saddle point formulation \eqref{conservative_prob} can be written as
\begin{align}
	\min_{\x \in \mathcal{X}} \max_{\lam \in \Rn_+^N} \ccalL(\x,\lam):= F(\x) + \ip{\lam,\H(\x) +\upsilon\textbf{1}}
	\label{lagrangian}
\end{align}
In this work, we build upon the Arrow-Hurwicz version of the saddle point algorithm. We consider the stochastic augmented Lagrangian $\L(\x,\lam,\th)$ given by  
\begin{align} \label{stochastic_lagrangian}
	\!\!\L(\x,\lam,\th)= f(\x,\th)+
	\sum_{i=1}^{N}\!\!\lambda_i(h_i\left(\x,\th\right)+\upsilon)
	-\frac{\delta \eta}{2} \norm{\lam}^2.
\end{align}
Note that the expression in \eqref{stochastic_lagrangian} is different from the standard Lagrangian due to the presence of the term $-\frac{\delta \eta}{2} \norm{\lam}^2 $. The augmentation allows us to control the value of $\lam$ and prevent it from growing too large through the use of parameter $\delta$. Indeed, the use of the augmented term allows us to bound the norm of the dual iterates, thus obviating the need for an explicit constraint on the domain of $\lam$. The compactness assumption on the domain of $\lam$ is a key requirement in classical stochastic saddle point algorithms \cite{nemirovski2005efficient}. Further details regarding the choice of the parameters $\eta$, $\delta$, and $\upsilon$ will be provided later.
\begin{algorithm}[t]
	\caption{Conservative Stochastic Optimization Algorithm}
	\label{alg:CSOA}
	\begin{algorithmic}
		\State {\bfseries Initialization:} $\x_1$, $\eta$, $\delta$, $\upsilon$ and $\lam_1=0$.
		\For{$t=1$ {\bfseries to} $T$}
		\State Update $\x_{t+1}$ and $\lam_{t+1}$ as 
		\begin{align*}
			\hspace{-1cm}&\x_{t+1} 
			= \px\Big[\x_t  - \eta  \nabla f (\x_t,\th_t) -\eta \sum_{i=1}^{N}\lambda_{i,t}[ \nx h_i  \left(\x_t, \th_t\right)]\Big] \\
			\hspace{-8mm}&\lambda_{i,t+1} = \Big[ (1 - \eta^2 \delta) \lambda_{i,t} + \eta \left(h_i\left(\x_{t},\th_{t}\right)+\upsilon\right)\Big]_+ 
		\end{align*}
		for all $i\in \{1, \ldots, N\}$.
		\EndFor
	\end{algorithmic}
\end{algorithm}
\subsection{Conservative Stochastic Optimization Algorithm (\textbf{CSOA})} 
The proposed \textbf{CSOA} algorithm is essentially an application of the saddle point algorithm to $\L$ in \eqref{stochastic_lagrangian}. Starting with an arbitrary $\x_1$ and $\lam_1 = 0$, the updates take the form
\begin{align} \label{iterate_primal}
	\x_{t+1}   &=\px\Big[ \x_t - \eta \nx \L (\x_t, \lam_t,\th_t)\Big] \;,  \\
	\lam_{t+1} &= \Big[\lam_t + \eta \nabla_{\lam} \L (\x_{t}, \lam_t, \th_t)  \Big]_+\;,
	\label{iterate_dual}
\end{align}
where $\nx\L (\x_t, \lam_t,\th_t)$ and $\nabla_{\lam} \L (\x_{t}, \lam_t, \th_t)$, are the primal and dual stochastic gradients of the augmented Lagrangian, with respect to $\x$ and $\lam$, respectively. Here, $\px(\x)$ represents the projection of the vector $\x$ on to the compact set $\mathcal{X}$ while $[\cdot]_+$ denotes the projection on to $\Rn^m_{+}$. The full algorithm is summarized in Algorithm \ref{alg:CSOA}.
\begin{algorithm}[t]
	\caption{Frank Wolfe Conservative Stochastic Optimization Algorithm}
	\label{algo_2} 
	\begin{algorithmic}[1]
		\State {\bfseries Initialization:}  $\x_{0}$, $\x_{1}$, $\eta$, $\delta$, $\rho$, $\upsilon$, $\d_{1}=\textbf{0}$ and $\bblam_{1}=\textbf{0}$.
		\For {$t=1$ {\bfseries to} $T$}
		\State Update gradient estimate  as
		\begin{align*}
			\d_t&=(1-\rho)\d_{t-1} +\nx\L(\x_{t},\lam_{t},\th_{t})- (1-\rho) \nx\L(\x_{t-1},\lam_{t-1},\th_{t})
		\end{align*}
		\State Calculate $\bbs_t$ as
		$\bbs_t=\argmin_{\bbs \in \mathcal{X}} \ip{\bbs,\d_t}$
		\State Update $\x_{t+1}$  as
		$\x_{t+1} 
		=\x_t +\eta(\bbs_t-\x_t)$
		\State Update dual variable for each $i\in \{1,2,\cdots,N\}$ as:			
		\begin{align*} 
			\lambda_{i,t+1} = \Big[ (1 - \eta^2 \delta) \lambda_{i,t} + \eta \left(h_i\left(\x_{t},\th_{t}\right)+\upsilon\right)\Big]_+  
		\end{align*}			
		\EndFor
	\end{algorithmic}
\end{algorithm}


\subsection{Frank-Wolfe Conservative Stochastic Optimization Algorithm (\textbf{FW-CSOA})}
The Frank-Wolfe version avoids the projection in \eqref{iterate_primal} and is well-suited to cases when such an operation is computationally expensive \cite{nemirovski2009robust}. In deterministic settings where the goal is to minimize $F(\x)$ over $\x \in \cX$, the FW updates entail finding a direction $\s_t$ that satisfies $\s_t = \arg\min_{\s \in \cX} \ip{\s,\nabla F(\x_t)}$. However, similar updates can not be utilized in  stochastic settings; given an unbiased stochastic gradient $\nabla f(\x_t,\th_t)$, the stochastic FW update along $\hat{\s}_t = \arg\min_{\s \in \cX} \ip{\s,\nabla f(\x_t,\th_t)}$ does not correspond to an FW update in expectation. It is instead necessary to obtain a (possibly biased) estimator $\d_t$ of gradient $\nabla F(\x_t)$ that has a lower variance as compared to $\nabla f(\x_t,\th_t)$. Such an estimator can be obtained either by employing a mini-batch \cite{hazan2012projection} or tracking the gradient using the recursion $\d_t = (1-\rho)\d_{t-1} + \rho\nabla F(\x_t,\th_t)$ \cite{mokhtari2018stochastic}. 

In order to adopt the FW framework to CSOA and obtain the desired centralized rates, we require a better tracking algorithm. To this end, we consider the recursive gradient tracking approach of \cite{cutkosky2019momentum}, which can be used to track the gradient of the augmented Lagrangian, and takes the form:
\begin{align}
	\d_t&=(1-\rho)\d_{t-1} +\nx\L(\x_{t},\lam_{t},\th_{t}) - (1-\rho) \nx\L(\x_{t-1},\lam_{t-1},\th_{t}). \label{nh}
\end{align}
Observe here that such a tracking process requires two stochastic gradients for each function at every iteration, namely, $\nabla f(\x_t,\th_t)$ and $\nabla f(\x_{t-1},\th_t)$. The gradient estimate $\d_t$ is subsequently utilized to obtain the direction $\bbs_t=\argmin_{\bbs \in \mathcal{X}} \ip{\bbs,\d_t}$ which is then used to update the primal variable as
\begin{align}\label{iterate_primal_FW}
	\x_{t+1} =\x_t +\eta(\bbs_t-\x_t).
\end{align}
With the dual update remaining the same as in \eqref{iterate_dual}, the full algorithm takes the form shown in  Algorithm \eqref{algo_2}. 
\section{Convergence Analysis}\label{sec:convergence}
This section studies the convergence rate results for the proposed algorithms. We establish the bound on the objective function optimality gap $F(\x_t)-F(\x^\star)$ and the constraint violation, both in expectation. We prove that the primal iterates converge in expectation to the optimal $F(\x^\star)$ at a rate of $\mathcal{O}(T^{-1/2})$ with zero average constraint violation. Let ($\x^\star,\lam^\star$) and ($\xu^\star,\lam^\star_\upsilon$) be the primal-dual optimal pairs for the problems in \eqref{prob_form} and \eqref{new_prob}, respectively. 

To prove these convergence results, let us first discuss some preliminary definitions and results which are crucial to the analysis. Let $\mathcal{F}_t$ denote the filtration which collects  the algorithm history
$\{\th_u, \x_u, \lam_u\}_{u=1}^{t-1}$. The primal and dual stochastic gradient{s} are the unbiased estimate of the true primal and dual gradients,  respectively, which implies that 
$\mbE\big[\nx\L(\x_t, \lam_{t},\th_t)|\mathcal{F}_t\big]= \nx\ccalL(\x_t, \lam_{t})$ 
and 
$\mbE\big[\nabla_{\lam}\L(\x_t, \lam_{t},\th_t)|\mathcal{F}_t\big]= \nabla_{\lam}\ccalL(\x_t, \lam_{t}).$
Next, we present the technical assumptions required for the convergence analysis.
\begin{assumption}\label{boundedgrad}
	\normalfont
	(\emph{Bounded stochastic gradients}) The stochastic gradients of the objective and constraint functions have bounded second moments, i.e.,  $\mbE\|\nx f(\x,\th)\|^2\leq \sigma_{f}^2$ and $\mbE\|\nx h_i(\x,\th)\|^2\leq\sigma_{h}^2$. 
\end{assumption}
\begin{assumption}\label{boundedh}
	\normalfont
	(\emph{Bounded constraint function}): The constraint functions $h_i\left(\x,\th\right)$ have bounded second moments, i.e., $$\mbE[h_i\left(\x,\th\right)^2]\leq \sigma_{\lam}^2.$$
\end{assumption}
These assumptions are useful in bounding the gradient of the Lagrangian, presented in the following corollary.
\begin{corollary}\label{cor-bounds}
	Under Assumptions \ref{boundedgrad}-\ref{boundedh}, it holds that
	\begin{align}
		\mbE\norm{\nx\L(\x,\lam,\th)}^2 &\leq 2B^2(1+\|\lam\|^2) \label{grad_norm_sq_x_zero},\\
		\mbE\norm{\nabla_{\lam}\L(\x,\lam,\th)}^2 & \leq 4N\sigma_{\lam}^2+4N\upsilon^2 + 2\delta^2\eta^2 \norm{\lam}^2 \label{grad_norm_sq_lam_zero}
	\end{align}
	where $B := \max(\sigma_f,\sigma_h \sqrt{N})$.
\end{corollary}
\noindent \textbf{Proof:}
From the properties of the norm, the left-hand side of \eqref{grad_norm_sq_x_zero} can be bounded as
\begin{align}
\mbE\norm{\nabla_{\x}\L(\x,\lam,\th)}^2 &\leq 2\mbE\norm{\nabla_{\x} f(\x,\th)}^2 + 2\mbE\norm{\sum_{i=1}^{N}\lambda_i\nabla_{\x} h_i(\x,\th)}^2 \nonumber\\&\leq 2 \sigma_f^2+ 2N\sigma_h^2\|\lam\|^2 \leq 2B^2(1+\|\lam\|^2).
\end{align}
where the second inequality follows from Assumption \ref{boundedgrad} and $B = \max(\sigma_f,\sigma_h \sqrt{N})$. Along similar lines, 
\begin{align}
\mbE\norm{\nabla_{\lam}\L(\x,\lam,\th)}^2 &\leq 2 \mbE\norm{\h(\x,\th) + \upsilon \mathbf{1}}^2+2\delta^2\eta^2\|\lam\|^2 \nonumber\\
&\leq 4N\sigma_{\lam}^2+4N\upsilon^2+2\delta^2\eta^2 \norm{\lam}^2.
\end{align}
where we have used the bound in Assumption \ref{boundedh}. 

Observe here that in Corollary \ref{cor-bounds}, the right-hand sides are not constants but functions of the dual variable $\lam$. These bounds will subsequently be used to establish the required convergence rates. Corollary \ref{cor-bounds} is a departure from the analysis of other saddle point variants where it is explicitly assumed that the right-hand sides of \eqref{grad_norm_sq_x_zero}-\eqref{grad_norm_sq_lam_zero} are constants \cite{madavan2019subgradient}. However, in the present case, the dual regularization term involving $\norm{\lam}^2$ allows us to make use of the bounds in Corollary \ref{cor-bounds} directly. 
\begin{assumption}\label{compact}
	\normalfont (\emph{Compact domain})
	The convex set $\mathcal{X}$ is compact, i.e.,  $\|\x-\bby\|\leq D$ for all  $\x,\bby \in \mathcal{X}$ for all $(\x,\y)\in \mathcal{X}$. Further, $\cX \subset \text{relint}(\text{dom} f \bigcap \cap_i \text{dom} h_i)$. 
\end{assumption}
Note that since $F$ is convex and defined over a compact set $\cX$, it follows from Assumption \ref{compact} that $F$ and  $H_i \;\forall\; i$ are also Lipschitz continuous. Let $G_f$ and $G_h$ denote the Lipschitz constant of $F$ and $H_i$ respectively. That is, $\abs{F(\x)-F(\y)}\leq G_f\norm{\x-\y}$ and  $\abs{H_i(\x)-H_i(\y)}\leq G_h\norm{\x-\y}$ for all $(\x,\y)\in \mathcal{X}$ and for all $i\in \{1, \ldots, N\}$.

For the constrained setting at hand, we also need a constraint qualification. In particular, we consider a strong version of Slater's condition. 
\begin{assumption}\label{slater}
	\normalfont (\emph{Strict feasibility})
	There exits a strictly feasible solution $\tilde{\x} \in \cX$ to \eqref{prob_form} that satisfies $H_i(\tilde{\x})+\sigma\leq 0$ for all $1 \leq i \leq N$ and for some $\sigma>0$.
\end{assumption}
It is remarked that Slater's condition is often necessary when studying constrained convex optimization problems. In the present case, Assumption \ref{slater} turns out to be critical for establishing the following corollary. 
\begin{corollary}\label{coro1}
	Under Assumptions \ref{compact} and \ref{slater} and for $0\leq \upsilon\leq \sigma/2$, it holds that:
	\begin{align}
		|F(\x^\star_{\upsilon})-F(\x^\star)|\leq \frac{2G_fD}{\sigma}\upsilon =: C\upsilon.
		\label{col}
	\end{align}
\end{corollary}
In \eqref{col}, we have introduced the constant $C := \frac{2G_fD}{\sigma}$. The proof of Corollary \ref{coro1} can be found in  \cite[Lemma 2]{thomdapu2019optimal}, and is provided in  Appendix \ref{proof_corollary1} for convenience. One implication of  Corollary \ref{coro1} is that the gap between the optimal values attained by \eqref{prob_form} and \eqref{conservative_prob} is $\O(\upsilon)$. 

Finally, we state the smoothness assumption{,} which is only required for the analysis of FW-CSOA. 
\begin{assumption}\label{smooth}
	\normalfont (\emph{Smoothness})
	The instantaneous objective function $f\left(\x,\th\right)$ and the constraint functions $\{h_i\left(\x,\th\right)\}_{i=1}^N$ are smooth, i.e., $\|\nabla f(\x,\th)-\nabla f(\bby,\th)\|\leq L_f\|\x-\bby\|$ and $\|\nabla h_i(\x,\th)-\nabla h_i(\bby,\th)\|\leq L_h\|\x-\bby\|$ for all $(\x,\y)\in \mathcal{X}$.
\end{assumption}
%
\subsection{Convergence rate of \textbf{CSOA}}
The convergence rate of Algorithm \ref{alg:CSOA} will be established under Assumptions \ref{boundedgrad}-\ref{slater}. We begin with stating Lemma \ref{lemma1}, which bounds the Lagrangian difference $\L(\x_t,\lam,\th_t) -\L(\x,\lam_t,\th_t) $ in terms of the primal and dual iterates. 
\begin{lemma}\label{lemma1}
	Let $(\x_t, \lam_t )$ be the sequence of iterates generated by Algorithm \ref{alg:CSOA}. Then, under Assumptions \ref{boundedgrad}-\ref{compact}, we have the bound
	\begin{align}\label{eq:lemma1_strong}
		\L(\x_t,\lam,\th_t) -\L(\x,\lam_t,\th_t)   
		&
		\leq  \frac{1}{2\eta} \big( \|\x_{t}-\x\|^2 -\|\x_{t+1}-\x\|^2\big)+\frac{1}{2\eta}\big(\norm{\lam_{t}-\lam}^2 - \norm{\lam_{t+1}-\lam}^2 \big)  \nonumber
		\\
		&\quad +\frac{\eta}{2}\norm{ \nx\L(\x_t,\lam_t,\th_t)}^2+\frac{\eta}{2} \norm{\nabla_{\lam}\L(\x_t,\lam_t,\th_t)}^2.
	\end{align}
\end{lemma}
The proof of Lemma \ref{lemma1} follows directly from the updates in \eqref{iterate_primal} and \eqref{iterate_dual} and from the convex-concave nature of the augmented Lagrangian. 

\noindent \textbf{Proof:}
We divide proof into two parts. In the first part, we establish an upper bound on 
$\L \left(\x_t,\lam_t,\th_t\right) -\L(\x, \lam_t,\th_t)$, while in the second part, we establish a lower bound on $\L (\x_t,\lam_t,\th_t) -\L(\x_t, \lam,\th_t)$. The final bound required in Lemma \ref{lemma1} would follow from combining these two bounds. 

\textit{Bound on} $\L (\x_t,\lam_t,\th_t) -\L(\x, \lam_t,\th_t)$: For any $\x \in \cX$, it holds from \eqref{iterate_primal} that
\begin{align}
\norm{\x_{t+1}-\x}^2=&\norm{\px[\x_t-\eta\nx\L(\x_t,\lam_t,\th_t)]-\x}^2\nonumber\\
\leq &\norm{\x_t-\eta\nx\L (\x_t,\lam_t,\th_t)-\x}^2\label{expand_x0}\\
=&\norm{\x_t-\x}^2-2\eta\ip{\nx\L (\x_t,\lam_t,\th_t),\x_t-\x}+\eta^2\norm{\nx\L (\x_t,\lam_t,\th_t)}^2,\label{expand_x}
\end{align}
where \eqref{expand_x0} follows from the non-expansiveness property of the projection operator $\px[\cdot]$ while \eqref{expand_x} is obtained by expanding the square. Since the Lagrangian is convex in $\x$, it holds that $
\L(\x_t,\lam_t,\th_t)-\L(\x, \lam_t,\th_t)\leq \ip{\nx\L (\x_t,\lam_t,\th_t),\x_t- \x}$, which upon substituting into \eqref{expand_x}, and rearranging, yields 
\begin{align}\label{first half_lemma1}
&\L(\x_t,\lam_t,\th_t) -\L(\x, \lam_t,\th_t)\leq \frac{1}{2\eta} \left( \norm{\x_t-\x}^2\!  -\!  \norm{\x_{t+1} - \x}^2\right)\!+\!\frac{\eta}{2}\| \nx\L_t(\x_t,\lam_t,\th_t)\|^2.
\end{align} 

\textit{Bound on} $\L (\x_t,\lam_t,\th_t) -\L(\x_t, \lam,\th_t) $: Proceeding along similar lines, for any $\lam \geq 0$, it follows from \eqref{iterate_dual} that 

\begin{align}
\norm{\lam_{t+1}-\lam}^2&=\norm{[\lam_t+\eta\nabla_{\lam}\L(\x_t,\lam_t,\th_t)]_{+}-\lam}^2
\label{temp:lam}\\
&\leq \norm{\lam_t+\eta\nabla_{\lam}\L(\x_t,\lam_t,\th_t)-\lam}^2\nonumber\\
&=\norm{\lam_t-\lam}^2+2\eta\ip{\nabla_{\lam}\L(\x_t,\lam_t,\th_t),\lam_t-\lam} +\eta^2\norm{\nabla_{\lam}\L(\x_t,\lam_t,\th_t)}^2.\label{above}
\end{align} 
Since the Lagrangian is concave with respect to $\lam$, we have that $\L(\x_t,\lam_t,\th_t) -\L(\x_t, \lam,\th_t)\geq  \ip{\nabla_{\lam}\L (\x_t,\lam_t,\th_t),\lam_t-\lam}$, which upon substituting into \eqref{above} and rearranging, yields
\begin{align}\label{second half_lemma1}
&\L(\x_t,\lam_t,\th_t) -\L(\x_t, \lam,\th_t)\geq \frac{1}{2\eta}\left(\norm{\lam_{t+1}\!-\!\lam}^2  \!- \! \norm{\lam_t \!-\! \lam}^2 \right) \!  -\!\frac{\eta}{2} \norm{\nabla_{\lam}\L(\x_t,\lam_t,\th_t)}^2
\end{align} 
Subtracting \eqref{second half_lemma1} from \eqref{first half_lemma1}, we obtain the required bound. 

Next, in order to establish the convergence proof for \textbf{CSOA}, we present the  following key lemma.
\begin{lemma}\label{lemma2}
	Let $(\x_t, \lam_t )$ be the sequence generated by Algorithm \ref{alg:CSOA}. Then for the choices $\eta < \frac{1}{4B}$, $\delta = 4B^2$, and $\upsilon < \sigma_\lam$, we have the bound:
	\begin{align}\label{main_equation_expctd_convex1}
		&\sum_{t=1}^{T}\mbE[F(\x_t)-F(\xus)] - \left(2B^2 \eta T+\frac{1}{2\eta}\right) \norm{\lam}^2
		+\sum_{t=1}^{T}\mbE[\ip{\lam,\H\left(\x_t\right)+\upsilon\one }]\leq \frac{D^2}{2\eta}  + P\eta T.
	\end{align}
	where $P:=B^2+4N\sigma^2_{\lam}$.
\end{lemma}
The proof of Lemma \ref{lemma2} entails using the definition of the augmented Lagrangian in \eqref{eq:lemma1_strong} as well as simplifying the resulting inequality by appropriately choosing the parameters. 

\noindent \textbf{Proof:}
Recalling the definition of the augmented Lagrangian in \eqref{stochastic_lagrangian}, we can write
\begin{align} 
\hat{\cL}(\x_t,\lam,\th_t)\! -\! \hat{\cL}(\x,\lam_t,\th_t)\! =\! f(\x_t,\th_t)\! +\!\ip{\lam,\h\left(\x_t,\th_t\right)\!+\!\upsilon\one} 
\!-\! \tfrac{\delta\eta}{2}\norm{\lam}^2\!-\! f(\x,\th_t) - \ip{\lam_t,\h\left(\x,\th_t\right)\!+\!\upsilon\one }\! +\! \tfrac{\delta \eta}{2} \norm{\lam_t}^2 \label{diff_eq}
\end{align}
for any $\x \in \cX$. Substituting the bound in Lemma \ref{lemma1}, taking $\x = \xus$, and rearranging,  we obtain 
\begin{align}\label{main_equationT1}
f(\x_t,\th_t)&-	f(\xus,\th_t)+ \frac{\delta \eta}{2} [\norm{\lam_t}^2-\norm{\lam}^2]+\ip{\lam,\h\left(\x_t,\th_t\right)+\upsilon\one }-\ip{\lam_t,\h\left(\xus,\th_t\right)+\upsilon\one }\nonumber\\&
\leq\! \frac{\eta}{2}\| \nx\L(\x_t,\lam_t,\th_t)\|^2\!+\!\frac{1}{2\eta} \left( \norm{\x_t\!-\!\xus}^2  \!- \! \norm{\x_{t+1} \!-\! \xus}^2\right)\nonumber\\&\quad
+\frac{\eta}{2} \|\nabla_{\lam}\L(\x_t,\lam_t,\th_t)\|^2+\frac{1}{2\eta}\left(\norm{\lam_t\!-\!\lam}^2 \! - \! \norm{\lam_{t+1} \!-\! \lam}^2 \right).
\end{align}
Also recall that $\F_t$ is the sigma field generated by $\{\th_u, \x_u, \lam_u\}_{u\leq t-1}$, so that $\mbE[f(\x_t,\th_t)-	f(\xus,\th_t)|\F_t]=F(\x_t)-F(\xus)$ and $\mbE[\h(\xus,\th_t)|\F_t]=\H(\xus)$. Therefore, taking expectation in \eqref{main_equationT1}, we obtain
\begin{align}\label{expected1}
\mbE[F(\x_t)-F(\xus)]& + \frac{\delta \eta}{2} \mbE[\norm{\lam_t}^2-\norm{\lam}^2]+\mbE[\ip{\lam,\H(\x_t)+\upsilon\one }] \nonumber\\
& \leq \frac{\eta}{2}\mbE\norm{\nx\L(\x_t,\lam_t,\th_t)}^2+\frac{\eta}{2} \mbE\norm{\nabla_{\lam}\L(\x_t,\lam_t,\th_t)}^2 \nonumber\\
&\quad + \frac{1}{2\eta} \E{ \norm{\x_t-\xus}^2  -  \norm{\x_{t+1} - \xus}^2} +\frac{1}{2\eta}\E{\norm{\lam_t-\lam}^2  -  \norm{\lam_{t+1} - \lam}^2 }
\end{align}
where we dropped the last term on the left since $\lam_t \geq 0$ and $\H(\xus) + \upsilon\one \leq 0$. Substituting the gradient bounds in Corollary \ref{cor-bounds} and rearranging, we obtain
\begin{align}\label{main_equation_expctd4}
\mbE[ F(\x_t)-F(\xus)] &- \frac{\delta \eta}{2} \norm{\lam}^2+\E{\ip{\lam,\H(\x_t)+\upsilon\one }}\nonumber\\
&\leq\frac{\eta}{2}[2B^2+4N\sigma^2_{\lam} + 4N\upsilon^2\! + (2B^2\! +2\delta^2\eta^2\!-\delta)\mbE\norm{\lam_t}^2] \nonumber\\
&\quad+  \frac{1}{2\eta} \mbE\Big[ \norm{\x_t-\xus}^2  -  \norm{\x_{t+1} - \xus}^2\Big]+\frac{1}{2\eta} \mbE\Big[\norm{\lam_t-\lam}^2  -  \norm{\lam_{t+1} - \lam}^2 \Big]. 
\end{align}
It is easy to verify that if $\eta < \frac{1}{4B}$, then the choice $\delta = 4B^2$ ensures that $2B^2 +2\delta^2\eta^2 \leq \delta$ and hence the term multiplying $\E{\norm{\lam_t}^2}$ can be dropped. 
Summing \eqref{main_equation_expctd4} over  $t=1,2,\ldots,T$, and recalling that $\lam_1 = 0$, we obtain
\begin{align}\label{main_equation_expctd6}
&\sum_{t=1}^T\mbE[(F(\x_t)-F(\xus))+\ip{\lam, \H(\x_t)+\upsilon\one}] - \frac{\delta \eta T}{2} \norm{\lam}^2\leq \frac{1}{2\eta} \mbE( \norm{\x_{1}-\xus}^2 + \norm{\lam}^2 )+ \frac{P\eta T}{2}.
\end{align}
where we have used the fact that $\upsilon < \sigma_\lam$ and substituted $P:=2B^2+8N\sigma^2_{\lam}$. Finally, the compactness of $\cX$ implies that $\norm{\x_1-\xus}^2 \leq D^2$, which upon substituting into \eqref{main_equation_expctd6} yields the desired bound. 

The next step is to choose $\lam$ in \eqref{main_equation_expctd_convex1} such that the left-hand side is maximized, resulting in a tighter bound, and subsequently, the required convergence rate result, as stated in the following Theorem. 
\begin{thm}\label{theorem1}
	For \textbf{CSOA}, under Assumptions \ref{boundedgrad}-\ref{slater}, and for the choices $\delta = 4B^2$, $\eta = C_1T^{-\frac{1}{2}}$, and $\upsilon = C_2T^{-\frac{1}{2}}$, the optimality gap is bounded as
	\begin{align}\label{eq:theorem1}
		\frac{1}{T} \sum_{t=1}^{T}\E{F(\x_{t})}-F(\x^\star)\leq KT^{-1/2}
	\end{align}
	and the cumulative constraint violation is bounded by zero, i.e., $\sum_{t=1}^{T} \E {H_i(\x_t)} \leq 0$ for all $i\in\{1,2,\cdots,N\}$. 	Here, the constants $C_1$, $C_2$, and $K$ depend only on the problem parameters as follows 
	\begin{align}
		C_1=&\mathcal{O}\left(\frac{1+D^2+(G_f^2D^2/\sigma^2)}{B^2+4N\sigma^2_\lambda+(16B^2G_f^2D^2/\sigma^2)}\right)\nonumber
		\\
		C_2=&\mathcal{O}\left((1+\sigma_{\lambda}\sqrt{N})\left(1+D+\frac{D}{\sigma}+\frac{D^2}{\sigma^2}+\frac{D^3}{\sigma^3}\right)\right)\nonumber
		\\
		K=&\mathcal{O}\left(C_1\left(1+B^2+D^2+N\sigma_\lambda^2+\frac{G_f^2D^2}{\sigma^2}\right)\right)
	\end{align}
which are derived from the expressions in \eqref{kc1c2}. 
\end{thm}
The proof of Theorem \ref{theorem1} is divided into two parts: establishing a bound on the optimality gap and establishing a bound on the constraint violation. The first part is largely similar to \cite{mahdavi2012trading} but modified to the stochastic setting. The second part is entirely novel and allows us to obtain the zero average constraint violation result. The approach is also different from \cite{bedi2019asynchronous}; in particular, the idea of bounding each constraint violation separately is taken from \cite{madavan2019subgradient}. 

As stated earlier, the idea is to tighten \eqref{main_equation_expctd_convex1} by appropriately choosing $\lam$. Once bounds on the optimality gap and the constraint violation have been obtained, the parameters $\upsilon$ and $\eta$ are adjusted to ensure the desired rates. It is clear that the left-hand side of \eqref{main_equation_expctd_convex1} is a concave quadratic function of $\lam$ and is maximized over $\Rn_{+}$ for the choice 
\begin{align}
	\hat{\lam}=\left(\frac{1}{4B^2 \eta T+ 1/\eta}\right)\left[\sum_{t=1}^T(\E{ \H(\x_t)}+\upsilon  \one)\right]_{+}.
	\label{opt_grad}
\end{align}
Substituting $\lam=\hat{\lam}$  into  \eqref{main_equation_expctd_convex1}, we obtain
	\begin{align}
		\sum_{t=1}^T \mbE[F(\x_t)&-F(\xus)] + \frac{\left[\sum_{t=1}^T(\E{ \H(\x_t)}+\upsilon  \one)\right]_{+}^2}{2(4B^2 \eta T+ \tfrac{1}{\eta})}\leq \frac{D^2}{2\eta} + \frac{P\eta T}{2}. \label{main_equation_expctd8}
	\end{align}
	Note that the second term on the left hand side of \eqref{main_equation_expctd8} is always positive, and hence by replacing it with its lower bound of zero, we could rewrite \eqref{main_equation_expctd8} as
	\begin{align}
		\sum_{t=1}^T \mbE[F(\x_t)-F(\xus)]
		&\leq \frac{D^2}{2\eta} + \frac{P\eta T}{2}. \label{main_equation_expctd822}
	\end{align}
	Adding and subtracting the optimal value $F(\x^\star) $ to the left hand side of \eqref{main_equation_expctd822}, we obtain 
	\begin{align}
		\sum_{t=1}^T \mbE[F(\x_t)&-F(\xus)+F(\x^{\star})-F(\x^{\star})] \leq \frac{D^2}{2\eta} + \frac{P\eta T}{2}. \label{main_equation_expctd823}
	\end{align}
	After rearranging the terms in \eqref{main_equation_expctd823}, we can write
	\begin{align}
		\sum_{t=1}^T \mbE[F(\x_t)-F(\x^{\star})] \leq& \frac{D^2}{2\eta} + \frac{P\eta T}{2}+T(F(\x^\star_{\upsilon})-F(\x^\star))
		\nonumber
		\\
		\leq & \frac{D^2}{2\eta} + \frac{P\eta T}{2}+T(|F(\x^\star_{\upsilon})-F(\x^\star)|), \label{main_equation_expctd824}
	\end{align}
	where we have utilized the inequality $a\leq |a|$ for any $a\in\mathbb{R}$. Next, we utilize the upper bound $	|F(\x^\star_{\upsilon})-F(\x^\star)|\leq C\upsilon$ from Corollary \ref{coro1} into \eqref{main_equation_expctd824} to obtain
	\begin{align}\label{use_coro}
		\sum_{t=1}^T \E{F(\x_t)}-F(\x^\star) &\leq \frac{D^2}{2\eta} + \frac{P\eta T}{2} + C\upsilon T.
\end{align}
Now we turn to establish {the} bound on constraint violation.  Taking total expectation in \eqref{eq:lemma1_strong},  utilizing the bounds in Corollary \ref{cor-bounds}, and substituting $P=2B^2+8N\sigma^2_{\lam}$,  we obtain:
\begin{align}\label{main_equation_expctd444}
	&\mbE[ \L(\x_t,\lam) -\L(\x, \lam_t)]\nonumber\\&\leq \frac{1}{2\eta} \mbE\Big[ \norm{\x_t-\x}^2  -  \norm{\x_{t+1} - \x}^2+\norm{\lam_t-\lam}^2 -  \norm{\lam_{t+1} - \lam}^2 \Big] +\frac{\eta}{2}[P+2(B^2 +\delta^2\eta^2)\mbE\norm{\lam_t}^2].
\end{align}
for any $\x \in \cX$ and $\lam \geq 0$. Here we have defined the averaged version of the augmented Lagrangian as $\L(\x,\lam) = \Ex{\L(\x,\lam,\th)} = F(\x) + \ip{\lam,\H(\x)+\upsilon\one} - \frac{\delta\eta}{2}\norm{\lam}^2$ where the expectation is with respect to $\th$. Comparing the standard Lagrangian $\cL(\x_t,\lam) $ (defined in \eqref{lagrangian}) with the expected augmented Lagrangian $\L(\x,\lam)$, we have an extra term 
$-\frac{\delta \eta}{2} \norm{\lam}^2$ so that $\L(\x,\lam) = \cL(\x,\lam) -  \frac{\delta\eta}{2}\norm{\lam}^2$. Thus from \eqref{main_equation_expctd444} we can upper bound $\mbE[ \cL(\x_t,\lam) -\cL(\x, \lam_t)]$ as follows:  
\begin{align}\label{main_equation_expctd443}
	&\mbE[\cL(\x_t,\lam) -\cL(\x, \lam_t)] 
	\\
	\quad & \leq \frac{\delta \eta}{2} \norm{\lam}^2+
	\frac{1}{2\eta} \mbE\Big[\norm{\x_t-\x}^2  -  \norm{\x_{t+1} - \x}^2+\norm{\lam_t-\lam}^2-  \norm{\lam_{t+1} - \lam}^2 \Big]+\frac{\eta}{2}[P+(2B^2 +2\delta^2\eta^2-\delta)\mbE\norm{\lam_t}^2].\nonumber
\end{align}
Let us first concentrate on the term $(2B^2 +2\delta^2\eta^2-\delta)$ which is quadratic in $\delta$.
	We note that for the selection $\delta = 4B^2$ and $\eta < \frac{1}{4B}$, we have that  $2B^2 +2\delta^2\eta^2-\delta<0$, which implies that the quadratic term multiplying $\E{\norm{\lam_t}^2}$ can be dropped from the right hand side of \eqref{main_equation_expctd443}.  Hence the inequality in \eqref{main_equation_expctd443} simplifies to 
	\begin{align}\label{main_equation_expctd4431}
		&\mbE[\cL(\x_t,\lam) -\cL(\x, \lam_t)] 
		 \leq {(2 \eta B^2)} \norm{\lam}^2+
		\frac{1}{2\eta} \mbE\Big[\norm{\x_t-\x}^2  -  \norm{\x_{t+1} - \x}^2+\norm{\lam_t-\lam}^2-  \norm{\lam_{t+1} - \lam}^2 \Big]+\frac{\eta P}{2}.
	\end{align}
	Taking summation on both sides of \eqref{main_equation_expctd4431} from $t=1,2,\cdots,T$, we obtain 
	\begin{align}\label{main_equation_expctd4432}
		&\sum_{t=1}^T\mbE[\cL(\x_t,\lam) -\cL(\x, \lam_t)] 
	 \nonumber\\& \leq {(2 \eta B^2 T)} \norm{\lam}^2+
		\frac{1}{2\eta} \mbE\Big[\norm{\x_1-\x}^2  -  \norm{\x_{T+1} - \x}^2+\norm{\lam_1-\lam}^2-  \norm{\lam_{T+1} - \lam}^2 \Big]+\frac{\eta P T}{2},
	\end{align}
	where the terms including primal and dual variables on the right of \eqref{main_equation_expctd4432} cancel out due to the telescopic sum. Next, we replace the negative terms on the RHS of \eqref{main_equation_expctd4432} with their upper bound of $0$ and substitute  $\lam_1=\mathbf{0}$ to obtain
	\begin{align}\label{main_equation_expctd4433}
		&\sum_{t=1}^T\mbE[\cL(\x_t,\lam) -\cL(\x, \lam_t)] 
	 \leq {(2 \eta B^2T)} \norm{\lam}^2+	\frac{\norm{\lam}^2}{2\eta}+
		\frac{1}{2\eta} \mbE\big[\norm{\x_1-\x}^2\big]+\frac{\eta P T}{2}.
	\end{align}
	It remains to substitute $\x=\xus$ and utilize the bound $\|\x_1-\xus\|\leq D$ into \eqref{main_equation_expctd4433}, which yields
	\begin{align}\label{lag_diffrnc_bound222}
		&\sum_{t=1}^T\mbE[\cL(\x_t,\lam) -\cL(\xus, \lam_t)] 
	 \leq 
		\frac{D^2}{2\eta}+\left(2 \eta B^2T+\frac{1}{2\eta}\right)\norm{\lam}^2+\frac{\eta P T}{2}.
\end{align}
Let $(\xus,\lus)$ be the primal-dual optimal pair for the problem \eqref{new_prob} and let $\one_i$ be such that $[\one_i]_j = 1$ for $i = j$ and zero otherwise. Then it can be seen that
\begin{align} 
	&\E{\cL(\x_t,\one _i+\lus)} =\E{F(\x_t)}+
	\ip{\one _i+\lus, \E{\H(\x_t)}+\upsilon\one }\nonumber\\
	&=\E{F(\x_t)}+\ip{\lus,\E{\H(\x_t)}+\upsilon\one}+\E{H_i(\x_t)}+\upsilon\nonumber\\
	&=\E{\cL(\x_t,\lus)}+\E{H_i(\x_t)}+\upsilon. \label{stochastic_lagrangian223}
\end{align}
Since ($\xus, \lus$) is a saddle point of the convex-concave function $\cL(\x,\lam)$, it holds that  
\begin{align}\label{lag_bound_trick}
	\cL(\xus,\lam)\leq \cL(\xus,\lus) \leq \cL(\x,\lus)
\end{align}
for any $\x \in \cX$ and $\lam \geq 0$. Rearranging \eqref{stochastic_lagrangian223} and using \eqref{lag_bound_trick}, we obtain
\begin{align} \label{stochastic_lagrangian2233}
	\E{H_i(\x_t)}+\upsilon  &= \E{\cL(\x_t,\one _i+\lus) - \cL(\x_t,\lus)}\nonumber\\
	&\leq \E{\cL(\x_t,\one _i+\lus)-\cL(\xus,\lam_t)}.
\end{align}
Taking sum over $t=1,2,\ldots, T$ in \eqref{stochastic_lagrangian2233} and using the bound from \eqref{lag_diffrnc_bound222} with $\lam = \lus + \one_i$, we obtain
\begin{align} \label{stochastic_lagrangian553}
	\sum_{t=1}^T \E{H_i(\x_t)}+\upsilon T &\leq  \sum_{t=1}^T\E{\cL(\x_t,\one _i+\lus)-\cL(\xus,\lam_t)}\nonumber\\
	& \leq \tfrac{D^2}{2\eta} + \frac{P\eta T}{2} + \left(2B^2\eta T + \tfrac{1}{2\eta}\right)\norm{\lus+\one_i}^2
\end{align}
The last term in \eqref{stochastic_lagrangian553} can be bounded from \eqref{new_lag3}, which implies that
\begin{align}
	\norm{\lus+\one_i}^2  &= 1 + \norm{\lus}^2 \leq 1 + \norm{\lus}_1^2 \nonumber\\
	&\leq 1 +( \one^\top \lus)^2 \leq 1 + C^2.
\end{align}
Summarizing, we have the following bounds on the optimality gap and the cumulative constraint violation
\begin{align}
	&\sum_{t=1}^T \E{F(\x_t)}-F(\x^\star) \leq \frac{D^2}{2\eta} + \frac{P\eta T}{2} + C\upsilon T \\
	&\sum_{t=1}^T \E{H_i(\x_t)}+\upsilon T \leq \tfrac{D^2}{2\eta} + \frac{P\eta T}{2} + \left(2B^2\eta T + \tfrac{1}{2\eta}\right)\left( 1 + C^2\right)
\end{align}
Therefore, in order to ensure that the constraint violation is less than or equal to zero, we choose 
\begin{align}
	\upsilon = \tfrac{D^2}{2\eta T} + \frac{P\eta}{2}  + \left(2B^2\eta  + \tfrac{1}{2\eta T}\right)\left( 1 + C^2\right)
\end{align}
which implies that the optimality gap is bounded by 
\begin{align}
	&\sum_{t=1}^T \E{F(\x_t)}-F(\x^\star) \leq   \frac{D^2}{2\eta} + \frac{P\eta T}{2}+ C\left(\tfrac{D^2}{2\eta } + \frac{P\eta T}{2}  + \left(2B^2\eta T + \tfrac{1}{2\eta }\right)\left( 1 + C^2\right)\right)= \frac{K_1}{2\eta} + \frac{K_2\eta T}{2} \label{gap1}
\end{align}
where $K_1 := (D^2 +1 +C^2)(1+C)$ and $K_2 := (P + 4B^2(1+C^2))(1+C)$. The bound in \eqref{gap1} can be minimized by choosing $\eta = \sqrt{K_1/K_2 T}$, which translates to
\begin{align}
	\frac{1}T\sum_{t=1}^T \E{F(\x_t)}-F(\x^\star) \leq \sqrt{\frac{K_1K_2}T}.
\end{align}
Equivalently, we denote 
\begin{align}\label{kc1c2}
	K\! =\! \sqrt{K_1K_2};\;  
	C_1\! =\! \sqrt{\frac{K_1}{K_2}}\; \text{and}\;
	C_2\! =\! \frac{K_1+K_2}{2(1+C)}\sqrt{\frac{K_1}{K_2}}
\end{align}
Remarkably, the optimality gap decays as $\O(T^{-\frac{1}{2}})$, which is the optimal rate even for unconstrained stochastic optimization problems with general convex objectives \cite{shamir2013stochastic}. Additionally, since the constraint violation is zero, the convergence rate obtained in Theorem \ref{theorem1} is therefore optimal in terms of its dependence on $T$. It must be noted that $\eta = \upsilon = \O(T^{-\frac{1}{2}})$ and hence satisfy the conditions required in Lemma \ref{lemma2} for $T$ sufficiently large. 
\subsection{Convergence rate of the \textbf{FW-CSOA}}
In this section, we develop the convergence rate for Algorithm \ref{algo_2} using Assumptions \ref{boundedgrad}-\ref{smooth}. The key to the analysis lies in characterizing the tracking property of $\d_t$ as explicated in the following lemma. 
\begin{lemma}\label{lemma3}
	Under Assumptions \ref{boundedgrad}-\ref{smooth}, the iterates generated by Algorithm  \ref{algo_2} satisfy
	\begin{align}
		\mbE\norm{\d_t - \nx\L(\x_{t},\lam_{t})}^2&\leq 8\rho^2 B^2(1+\mbE\norm{\lam_t}^2)+(1-\rho)^2\mbE \|\d_{t-1}-\nx\L(\x_{t-1},\lam_{t-1})\|^2\nonumber\\
		&\quad +2(1-\rho)^2(\eta^2D^2\mbE[L_{\lam_t}^2]+8\eta^2N^2\sigma_{h}^2\sigma_\lambda^2)
		\label{new_lemma1}
	\end{align}
	where  $L_{\lam_t}=L_f + L_h \norm{\lam_t}_1$. 
\end{lemma}
The proof of Lemma \ref{lemma3} is given in Appendix \ref{lemma3_proof} and follows along similar lines as Lemma 1 of \cite{cutkosky2019momentum}. A key difference here is that the bound on the right also depends on the norm of the dual iterate $\lam_t$, and is therefore not necessarily bounded. The presence of the dual iterate on the right makes the subsequent analysis very different. We first establish an upper bound on $\sum_{t=1}^{T}\mbE[\L(\x_t,\lam,\th_t)-\L(\x,\lam_t,\th_t)]$.
\begin{lemma}\label{lemma4}
	For \textbf{FW-CSOA}, then for the choice $\upsilon < \sigma_\lam$ and under Assumptions \ref{boundedgrad}-\ref{smooth}, it holds that
	\begin{align}\label{bound_frst_lemma}
		\sum_{t=1}^{T}&\mbE[\L(\x_t,\lam)-\L(\x,\lam_t)]
	\leq A_1+\frac{1}{2\eta}\norm{\lam}^2+\frac{\eta}{2}\left(\frac{2A_2}{\eta}+2\delta^2\eta^2\right)\sum_{t=1}^T \mbE\norm{\bblam_t}^2
	\end{align}
\end{lemma}
where, $A_1$ and $A_2$ are constants defined in \eqref{A_1} and \eqref{A_2}, respectively. The proof of Lemma \ref{lemma4} is given in Appendix \ref{lemma4_proof}.  We will use this property to show the convergence of primal iterates (Algorithm \ref{algo_2}) presented in the following theorem. 
\begin{thm}\label{theorem2}
	For \textbf{FW-CSOA}, under Assumptions \ref{boundedgrad}-\ref{smooth}, and for the choices $\delta=18L^2D^2$, $\eta=\hat{C_1}T^{-3/4}$  and $\upsilon =\hat{C_2}T^{-1/4}$, the optimality gap is bounded as
	\begin{align}\label{eq:theorem2}
		\frac{1}{T} \sum_{t=1}^{T}\E{F(\x_{t})}-F(\x^\star)\leq \hat{K} T^{-1/4}
	\end{align}
	while the cumulative constraint violation is bounded by zero, i.e., $\sum_{t=1}^{T} \E {H_i(\x_t)} \leq 0$ for all $i\in\{1,2,\cdots,N\}$. Here, $L=\max\{L_f,L_h\sqrt{N},1\}$ and the constants $\hat{C_1}$, $\hat{C_2}$, and $\hat{K}$ depend only on the problem parameters $G_f$, $G_h$, $D$, $B$, $\sigma_f$, $N$, and $\sigma_\lam$, and are defined in \eqref{final_cons}.
\end{thm}
The proof of Theorem \ref{theorem2} is given in Appendix \ref{proof_throrem2}. When solving problems with expectation constraints, the rate result in Theorem \ref{eq:theorem2} is the first of its kind among the class of projection-free algorithms. Nevertheless, we remark that  unlike in the projected case, the rate obtained here ($\O( T^{-1/4})$) is worse than the best projection-free rate achieved in the unconstrained setting ($\O( T^{-1/3})$). Whether it is possible to achieve $\O( T^{-1/3})$ rate in the context of projection-free stochastic constrained optimization remains an open question. Before concluding rate results, the following remark regarding the lower bound is due. 
\begin{remark}
	Observe that both Theorems \ref{theorem1} and \ref{theorem2} establish that the constraint violation is zero on average over a specified horizon. Since the constraint may be violated at some iterations, the objective function value may even be lower than $F(\x^\star)$. Thus, to complete the results, the following proposition establishes a lower bound on $F(\x_t)$, averaged over $1\leq t \leq T$.  
	\begin{proposition}\label{prop1} For $Q:=\left[\frac{(\hat{C}^2+4\delta)(C+r)^2+D^2\hat{C}^2+5\delta C^2}{2\hat{C}r}\right]$, we have that
		\begin{enumerate}[(a)]		
			\item  for \textbf{CSOA} with $\delta$,$\eta$ and $\upsilon$ as defined in Theorem \ref{theorem1} and $\hat{C}=C_1$, the optimality gap is lower bounded as
			\begin{align}
				\frac{1}{T}\sum_{t=1}^{T}\mbE[F(\x_t)-F(\x^\star)]\geq-C(Q+C_2)T^{-1/2}.
				\label{bound_12}
			\end{align}  
			\item for \textbf{FW-CSOA} with $\delta$,$\eta$ and $\upsilon$ as defined in Theorem \ref{theorem2} and $\hat{C}=\hat{C_1}$, the optimality gap is lower bounded as
			\begin{align}
				\frac{1}{T}\sum_{t=1}^{T}\mbE[F(\x_t)-F(\x^\star)]\geq& -C(Q+\hat{C_2})T^{-1/4},
				\label{bound_13}
			\end{align}  
		\end{enumerate}
	\end{proposition}
	The proof of proposition \ref{prop1} can be obtained along the lines as done in [\cite{nedic2009subgradient}, Proposition 2]. However, for the sake of completeness,  we provided the proof in the context of our problem in Appendix \ref{bound_dual_optimal}.
\end{remark}
\section{Numerical Experiments}\label{application} 
In this section, we compare the performance of the proposed conservative stochastic optimization algorithm (\textbf{CSOA}) and its projection-free version (\textbf{FW-CSOA}) with other  state-of-the-art algorithms.
\subsection{Fairness-constrained Classification }
We consider the problem of fair classification, where the training problem can be posed as that of finding a decision boundary that maximizes the average classification accuracy while ensuring that the predicted labels are uncorrelated with a pre-specified sensitive attribute. In particular, we consider the formulation from \cite{zafar2015fairness}, where the optimization problem is written as 
\begin{align}
	&\min_{\bbtheta} -\frac{1}{N}\sum_{i=1}^{N} (y_i\log(\sigma(\bbtheta^T\x_i)) + (1-y_i)\log(1-\sigma(\bbtheta^T\x_i)) \nonumber\\
	&\text{s.t.}\;\;\; -c \leq 
	\frac{1}{N}\sum_{i=1}^{N}(s_i-\bar{s})\bbtheta^T\x_i\leq c
	\label{app_prob_main}
\end{align}
where $\sigma(\cdot)$ denotes the sigmoid function, $\x_i$ denotes the feature vector corresponding to data point $i$, $y_i \in \{0,1\}$ denotes the corresponding label and $s_i$ denotes a sensitive attribute of interest. The formulation in \eqref{app_prob_main} restricts the correlation between the sensitive attribute and the predicted label $\bbtheta^T\x_i$ to be within the range $[-c,c]$.  In general, $c$ is required to be small as compared to the default correlation value obtained in the case of the standard (unconstrained) classifier. The classifiers such as the one in \eqref{app_prob_main} can be used for automated decision making tasks, such as hiring or promotion, by extracting the candidate's resume into a feature vector $\x_i$ and training the classifier using human-annotated historical data. As discussed in \cite{zafar2015fairness}, by choosing $c$ appropriately, it is possible to prevent the classifier from learning historical biases that might be inherent to the training data, such as the candidate's sensitive attributes (race or gender). The fairness of a decision making process is often characterized in terms of disparate impact, which indicates the bias of decision making process towards a particular class. Usually, $p\%$-rule is used to quantify disparate impact, which states that the ratio between the percentage of examples with a certain attribute value classified as positive in decision outcome and the percentage of examples without that attribute value also classified as positive in decision outcome must not be less
than $p\%$ \cite{zafar2015fairness}. As it may be challenging
to directly incorporate the $p\%$ rule into the problem formulation, we introduce a parameter $c$ within the constraints to control the value of $p$.  The choice of $c$ directly influences the disparate impact of the hiring process. For instance, $c$ can be chosen to ensure that the hiring process satisfy the so-called 80\% rule (or more generally, the $p$\% rule), wherein the selection rate for candidates with $s_i = 1$ be at least 80\% of the selection rate for candidates with $s_i = 0$. 
\begin{figure}[h]
	\centering
	\includegraphics[scale=0.4]
	{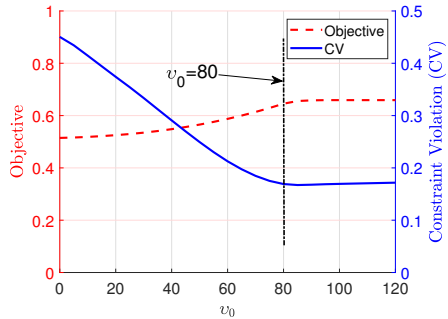}
	\label{cons_var_up} 
	\caption{The objective function value and the constraint violation on the Adult income dataset for different values of $\upsilon_0$, where we set $\upsilon = \upsilon_0 \times 1/\sqrt{T}$.}
	\vspace{-0mm}
	\label{fig:var_up}
\end{figure}

It can be observed that the problem in \eqref{app_prob_main} is an example of \eqref{prob_form}. However, the fair classifiers trained in \cite{zafar2015fairness} used batch algorithms, which are generally not scalable. We solve the problem in \eqref{app_prob_main} using the proposed \textbf{CSOA}  and compare its performance with state-of-the-art algorithm \textbf{G-OCO} \cite{yu2017online}. We consider the Adult income dataset \cite{asuncion20019uci}, which has a total of $ 45000$ examples, each with $14$ features such as education level, age, and marital status. The hiring process is simulated by attempting to predict if a subject has {an} income that exceeds $\$ 50,000$. We consider sex (male/female) as a sensitive attribute. In our experiments, we randomly sample $70\%$ of data to create training ($28350$ samples) and validation ($3150$ samples) sets, while the rest $30\%$ is used for testing the accuracy of the learned classifier. In order to set the problem parameter $c$, we randomly sample $10,000$ examples and solve the problem \eqref{app_prob_main} using the batch method \cite{zafar2015fairness}. We found that the choice $c = 0.05$ ensures that the $80\%$ rule is satisfied approximately and used this value throughout the simulations. When tuning the algorithm parameters as well, we ensure that the $p = 80\%$ is satisfied on the validation set. 

As dictated by theory, we set the parameters $\eta = \eta_0/\sqrt{T}$ and $\upsilon = \upsilon_0/\sqrt{T}$. Next, we tune the parameters $\eta_0$ and $\delta$ to achieve the best possible optimality gap on the training dataset, resulting in $\eta_0 = 0.15$ and $\delta = 10^{-2}$. Finally, we tune $\upsilon$ to minimize the constraint violation without significantly sacrificing the optimality gap. This was indeed possible, as evident from  Fig. \ref{fig:var_up}, where we set  $\upsilon = \upsilon_0\times /\sqrt{T}$ for different values of $\upsilon_0$. Finally, we selected $\upsilon_0 = 80$. The parameters for the \textbf{G-OCO} algorithm [cf. Algorithm 1 \cite{yu2017online}] were tuned in a similar manner, resulting in $\alpha=2.5\times T$ and $V=0.5\times \sqrt{T}$. 
\begin{figure*}[h]
	\centering
	\setcounter{subfigure}{0}
	\begin{subfigure}{0.33\columnwidth}		\includegraphics[width=\linewidth, height = 0.7\linewidth]{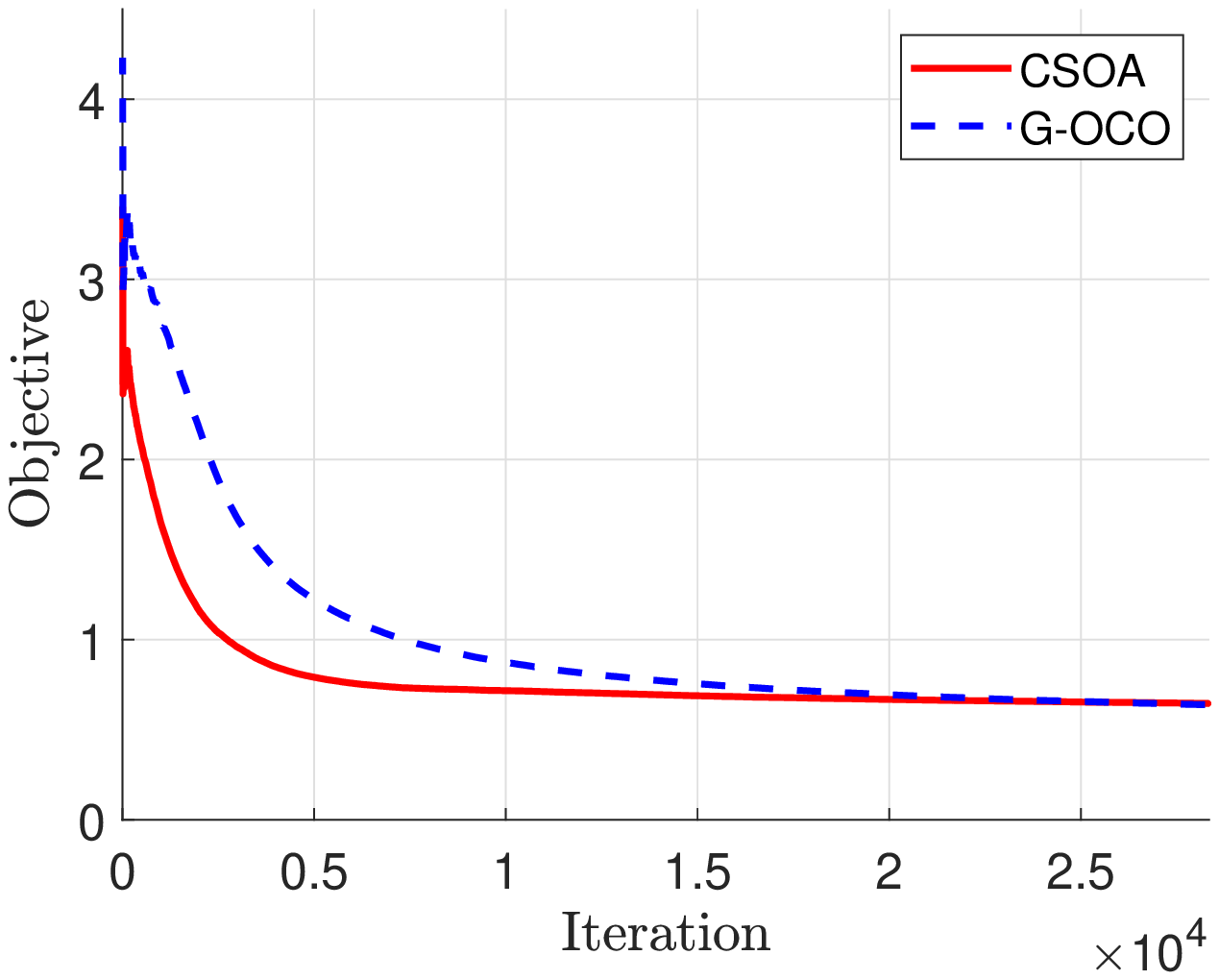}
		\caption{Objective function value}
		\label{NET1}
	\end{subfigure} 	
	\begin{subfigure}{0.33\columnwidth}
		\includegraphics[width=\linewidth, height = 0.7\linewidth]{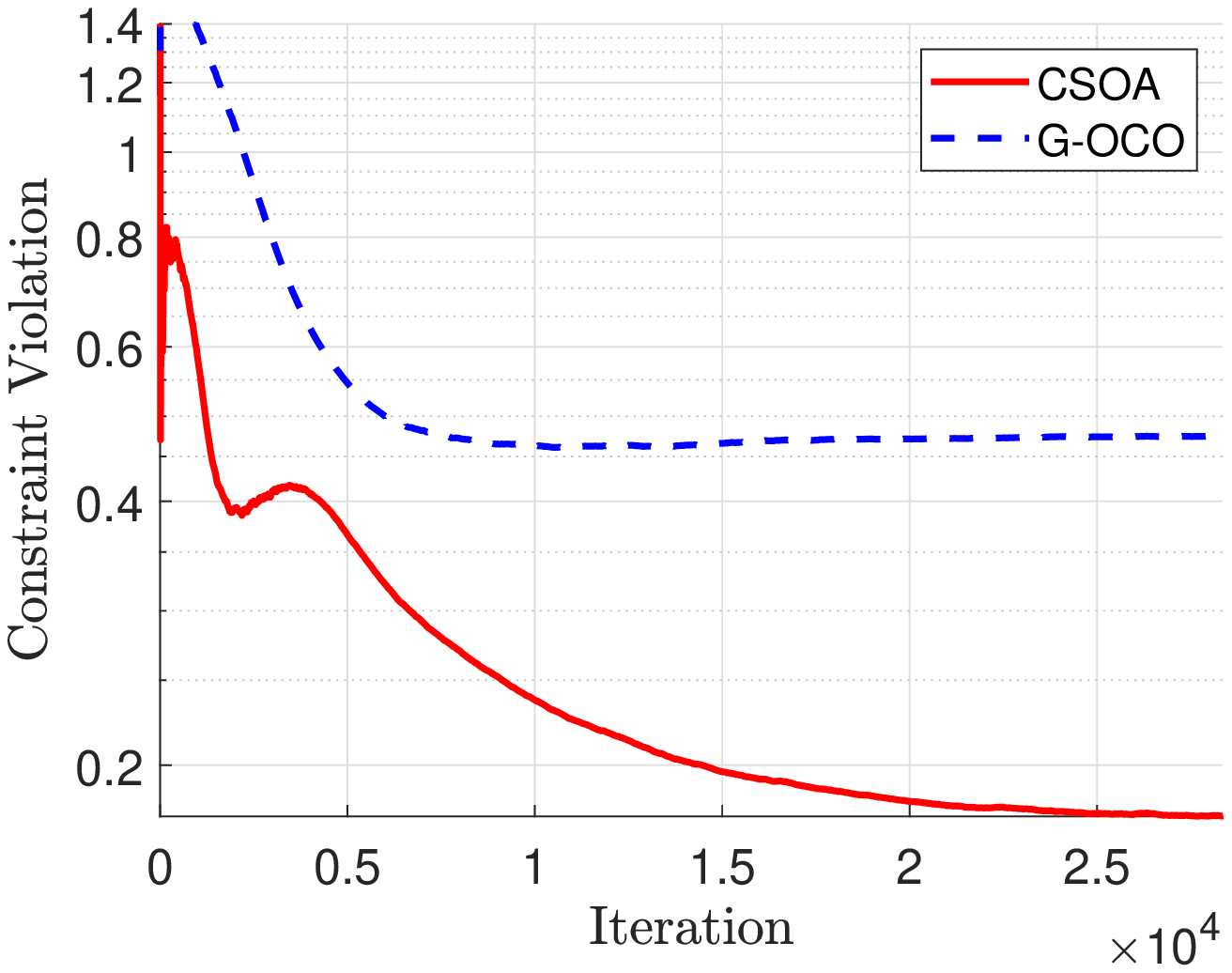}
		\caption{Constraint violation}
		\label{CVT1}
	\end{subfigure}
	\begin{subfigure}{0.33\columnwidth}
		\includegraphics[width=\linewidth, height = 0.7\linewidth]{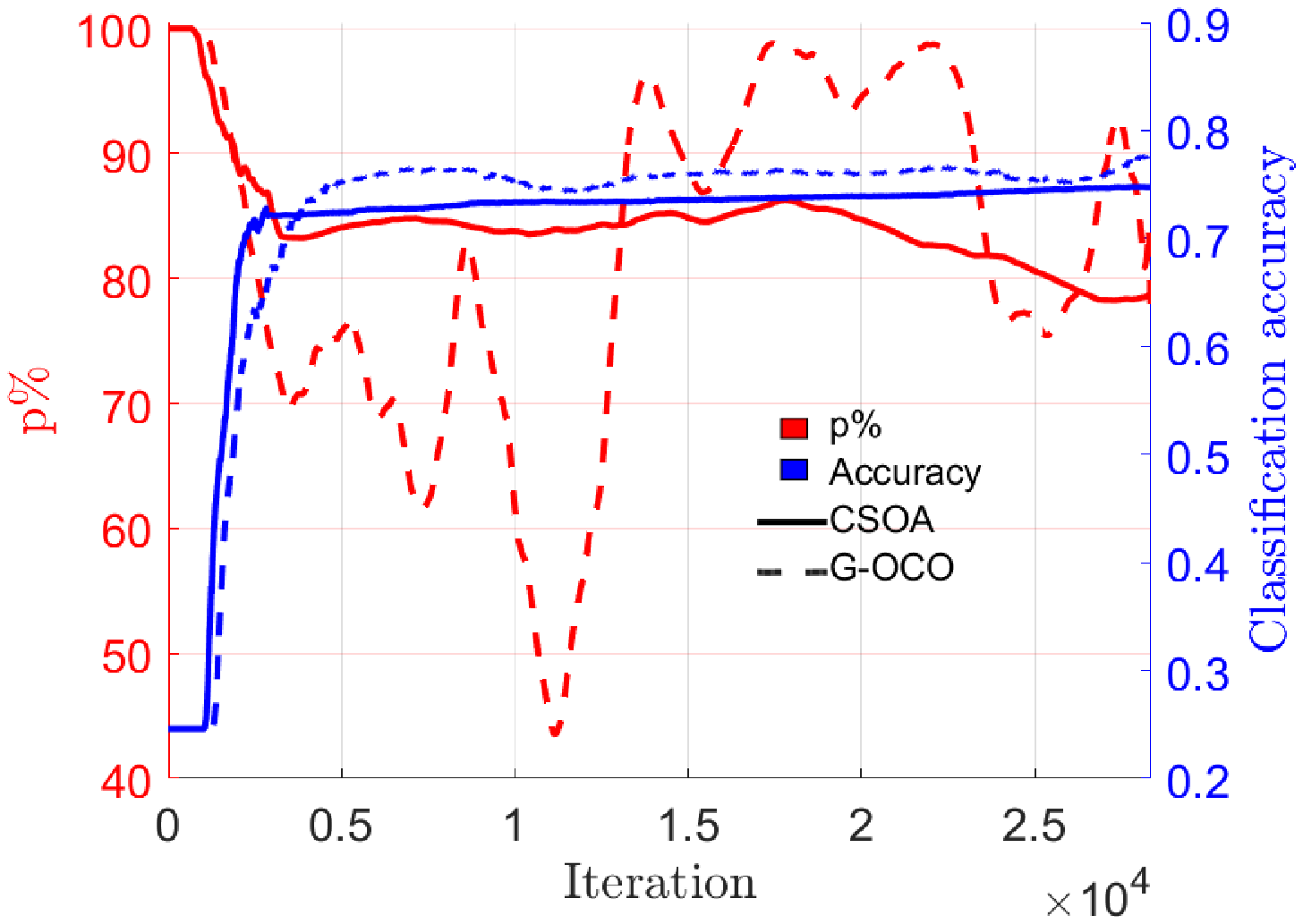}
		\caption{Accuracy and $p\%$ }
		\label{ACCT1}
	\end{subfigure}
	\vspace{0.25cm}
	\caption{Comparison of \textbf{G-OCO} (dotted lines) and the proposed \textbf{CSOA} method (solid lines). }
	\label{fig:fairness}
\end{figure*}

Fig.~\ref{fig:fairness} shows the performance of the two algorithms. It can be seen that while the objective function value evolves in a similar fashion for both the algorithms, the constraint violation for the proposed algorithm decays to almost zero. Such behavior is possible and expected due to the appropriate choice of the tightening parameter $\upsilon$. Fig.~\ref{ACCT1} compares the evolution of the classification accuracy (blue) and $p\%$ value (red) with iterations, both calculated on the test set. Observe here that the classification accuracy of the \textbf{CSOA} is almost the same as that of the \textbf{G-OCO}, despite the proposed algorithm satisfying the constraint more strictly. Further, it can be seen that the $p\%$ value of the \textbf{G-OCO} exhibits huge variations despite using exactly the same $c$ value as the proposed algorithm. In contrast, the proposed algorithm exhibits relatively stable behavior. 

We remark that as in \cite{zafar2015fairness}, it is again possible to study the trade-off between fairness and accuracy, and the same is reported in Appendix \ref{extra_expp}.
\begin{figure}[h!]
	\centering
	\setcounter{subfigure}{0}
	\begin{subfigure}{0.33\columnwidth}		\includegraphics[width=\linewidth, height = 0.7\linewidth]{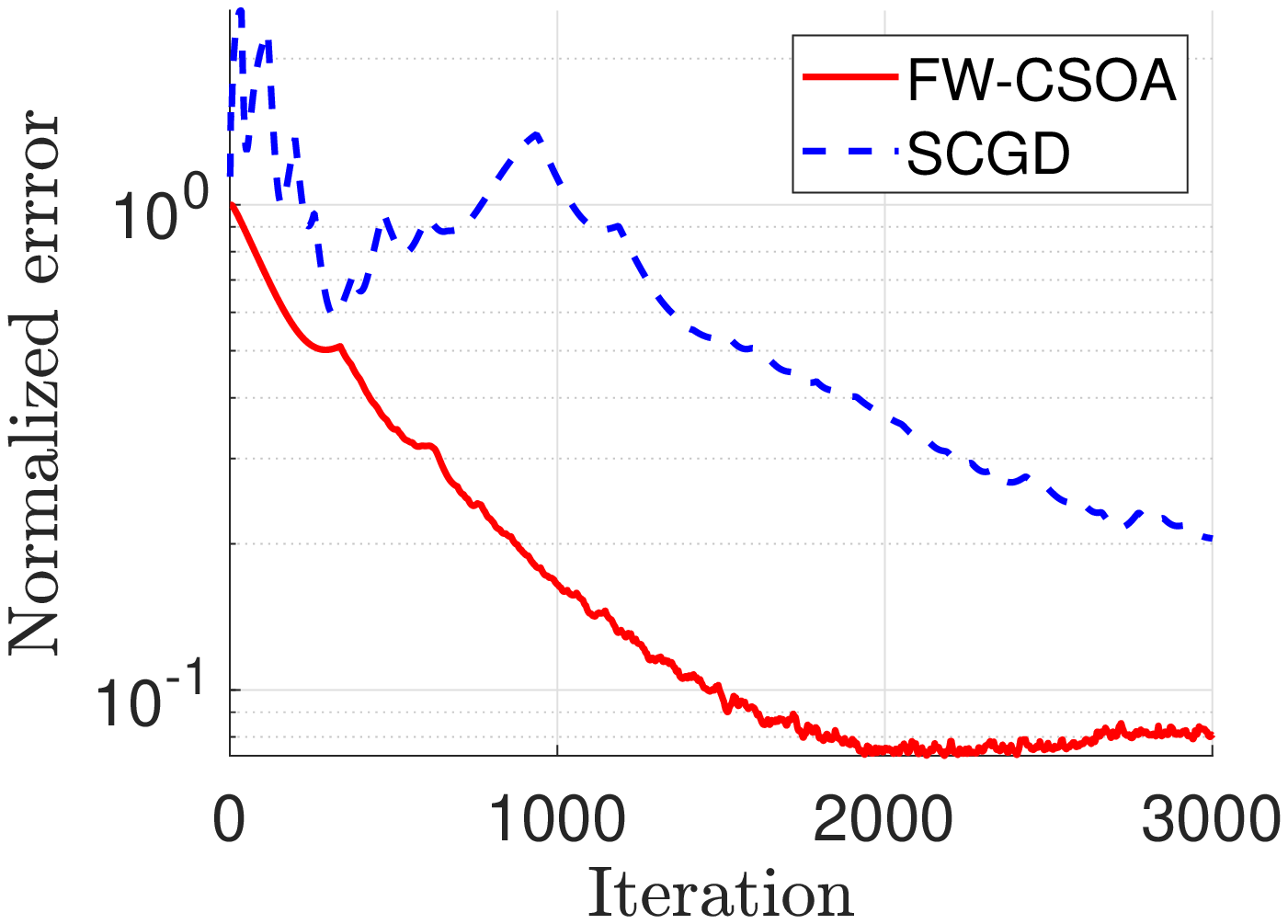}
		\caption{Objective function}
		\label{NET} 
	\end{subfigure}
	\begin{subfigure}{0.33\columnwidth}
		\includegraphics[width=\linewidth, height = 0.7\linewidth]{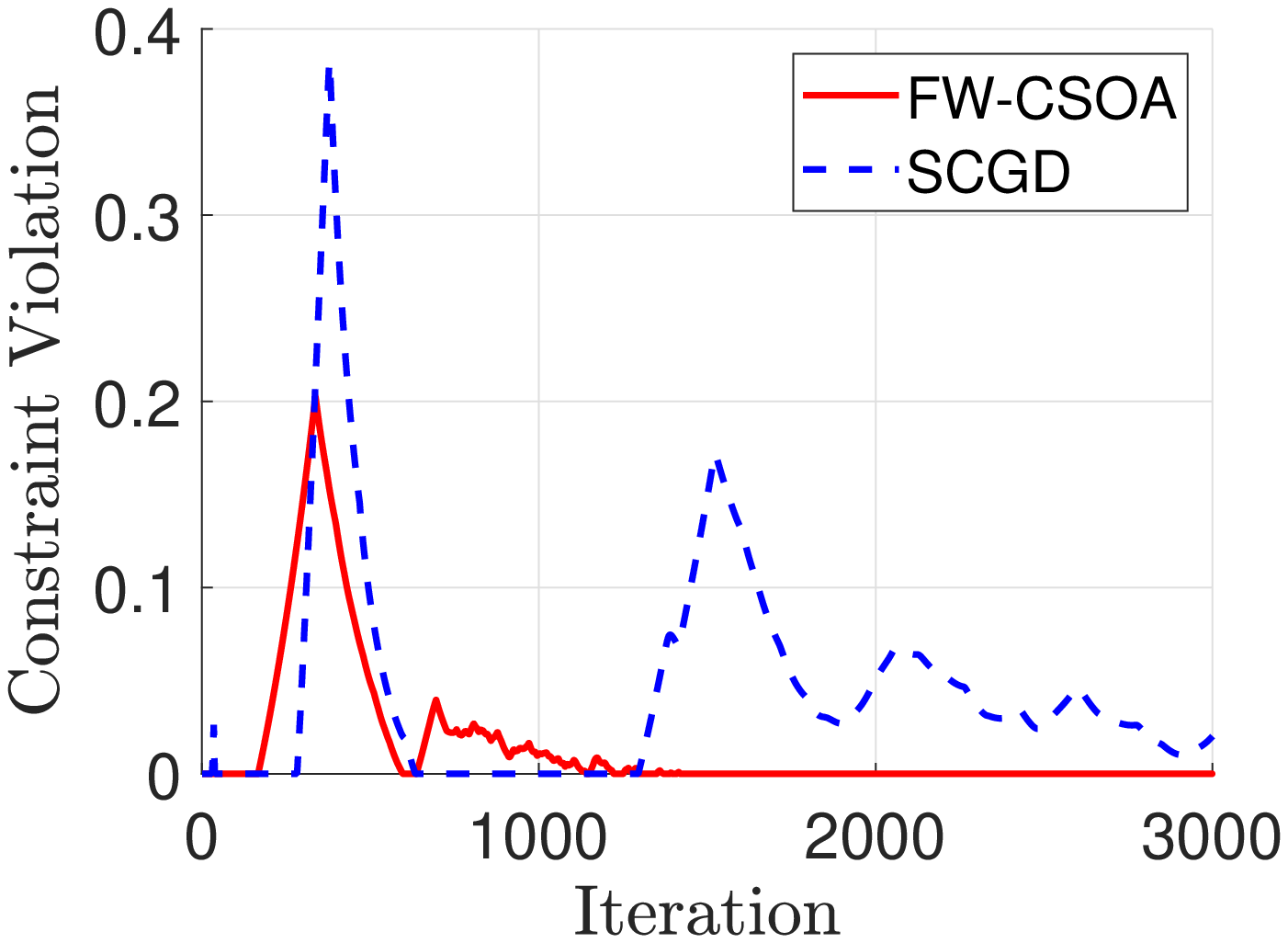}
		\caption{Constraint violation}
		\label{CVT} 
	\end{subfigure}
	\begin{subfigure}{0.33\columnwidth}		\includegraphics[width=\linewidth, height = 0.7\linewidth]{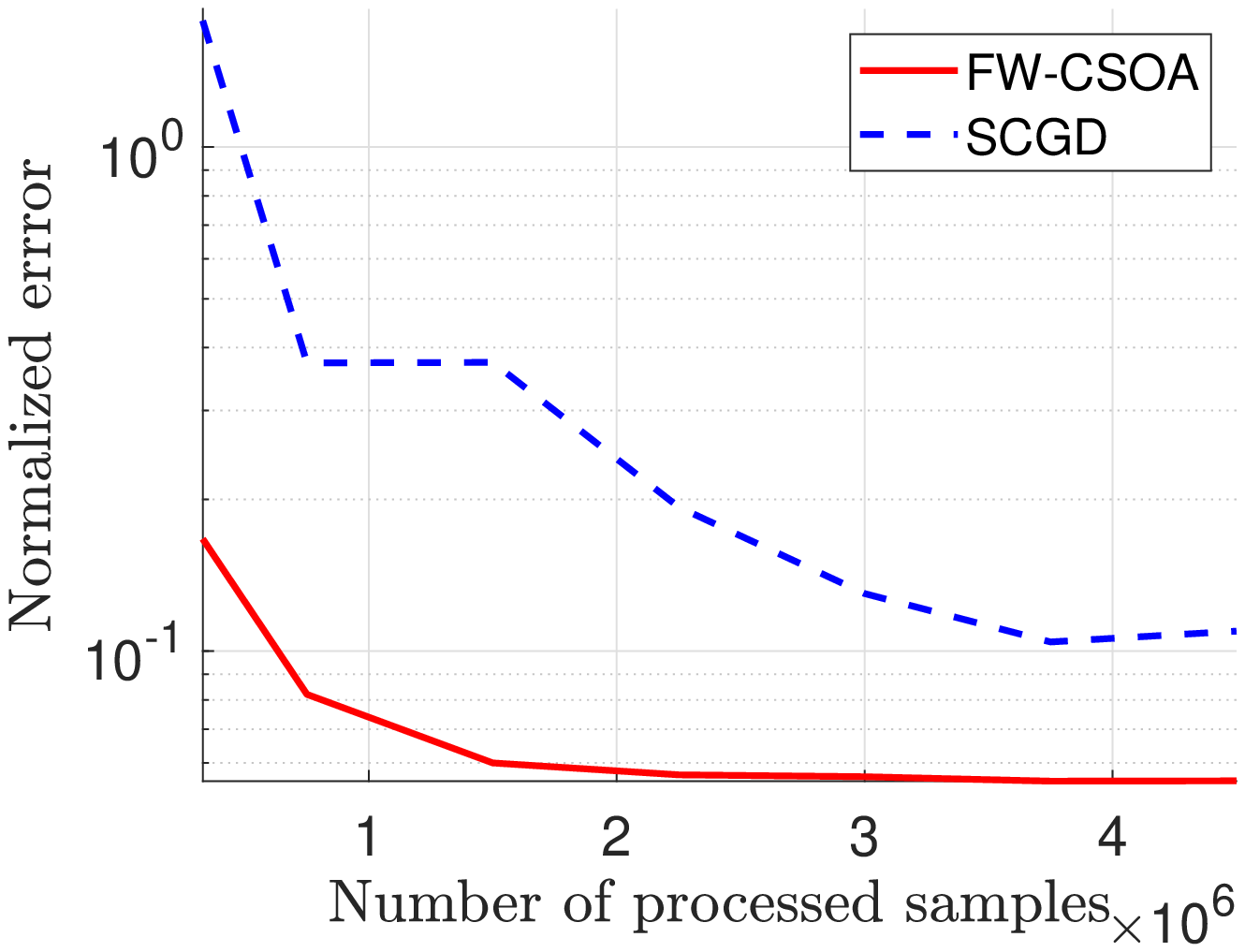}
		\caption{Objective function}
		\label{NES} 
	\end{subfigure}
	\begin{subfigure}{0.33\columnwidth}
		\includegraphics[width=\linewidth, height = 0.7\linewidth]{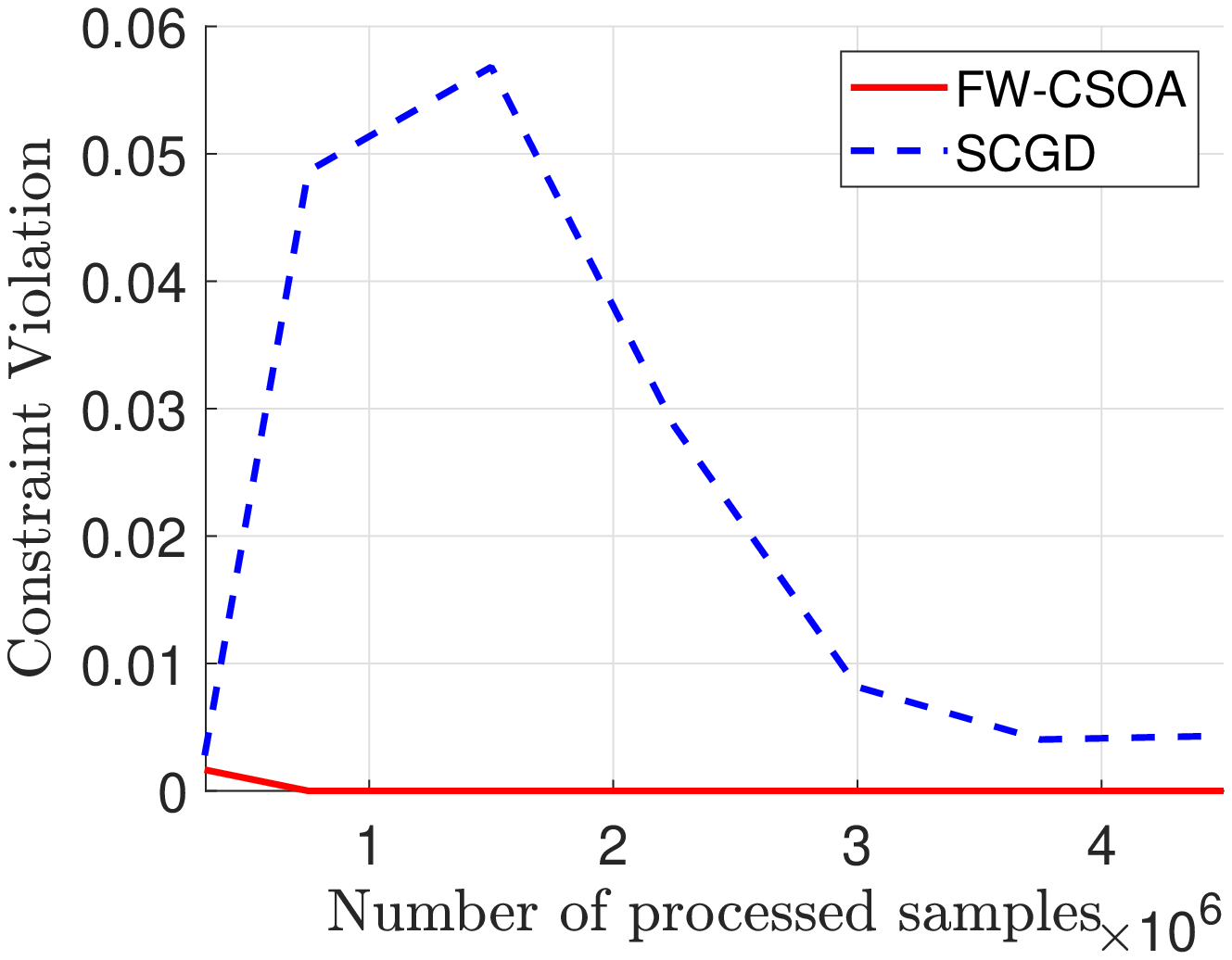}
		\caption{Constraint violation}
		\label{CVS} 
	\end{subfigure}
	\caption{Comparison of proposed algorithm FW-CSOA and SCGD \cite{mokhtari2018stochastic} for matrix completion problem discussed in  \ref{mat_comp}. Plot (a)-(b): Comparison in terms of number of iterations. Plot (c)-(d): Comparison  in terms of number of samples. Observe that FW-CSOA achieves much smaller error while maintaining very low constraint violation as compared to SCGD.}
	\vspace{-0mm}
	\label{fig:matrix_comp}
\end{figure}
\subsection{Matrix Completion for Structured Observations}\label{mat_comp}
We study the performance of the proposed projection-free \textbf{FW-CSOA} algorithm on the structural matrix completion problem. Matrix completion is a well-studied problem and has been immensely successful in the design of recommendation systems. Given a partially filled matrix, with entries denoting say user ratings for different movies, the goal is to complete the matrix, thereby predicting the ratings for movies not seen by the users. Traditionally, matrix completion is performed under the assumption that the values of the missing entries are independent of corresponding locations. In most recommender systems, however, such an assumption may not hold. In the context of predicting movie ratings, for instance, the fact that a user has not seen a particular movie may indicate the user's lack of interest in that movie. In structural matrix completion, we seek to exploit such a correlation in order to predict better the missing entries  \cite{molitor2018matrix}. 

Let $\M \in \Rn^{m \times n}$ be the observation matrix{,} and only a subset $\cI$ of the entries of $\M$ are observed while the entries in $\cI^c$ are missing. The goal is to find $\X$ that is as close as possible to $\M$ on the observed entries $(i,j)\in\mathcal{I}$. The missing entries are imputed by invoking the low-rank assumption on $\X$. To keep the problem convex, we simply use the nuclear norm relaxation and impose the constraint $\norm{\X}_\star \leq \alpha$ for some problem parameter $\alpha$. In order to account for such structural differences between observed and missing entries of the matrix, it was suggested in \cite{molitor2018matrix} that the values of the missing entries be explicitly forced to be small. Overall, the structural matrix completion can be posed as the following optimization problem:
\begin{subequations}
	\begin{align}\label{matrix_completion}
		&\min f(\bbX) : = \frac{1}{2} \sum_{i,j\in \cI} (X_{ij}-M_{ij})^2, 
		\\ & \text{s.t.} \hspace{10mm} \;\;\norm{\X}_{\star} \leq \alpha \label{CV1} 
		\\ &\hspace{8mm}  \;\; \frac{1}{2} \sum_{i,j\in \cI^c} X_{ij}^2 \leq \beta, 
		\label{CV2}
	\end{align}
\end{subequations}
where $X_{ij}$ represents the individual element $[\X]_{i,j}$ of matrix $\X$ at position $(i,j)$. 

For the simulations,  we follow the data generation process described in \cite{molitor2018matrix}. We consider the low rank matrix ${\X}^\star=\mathbf{X}_L^\star\mathbf{X}_R^\star$, where $\mathbf{X}_L^\star\in \mathbb{R}^{m\times r}$ and $\mathbf{X}_R^\star\in \mathbb{R}^{r\times n}$ are sparse matrices whose entries are chosen to be zero with probabilities $0.7$ and $0.5$, respectively. The non-zero entries are again uniformly distributed between 0 and 1. Set $\cI$ is obtained by subsampling $1\%$ of the zero and $90\%$ of the non-zero entries of ${\X}^\star$, and the observed entries are set as $[\M]_{i,j\in \cI}=[{\X}^\star]_{i,j\in \cI}+[\mathbf{W}]_{i,j\in \cI}$, where $\mathbf{W}$ is an error matrix generated as $[\mathbf{W}]_{i,j}=\gamma\cdot\frac{\norm{\X_{\star}}_F}{\norm{\mathbf{Z}}_F}\cdot[\mathbf{Z}]_{i,j}$, where $\gamma$ is a noise parameter and $\mathbf{Z}\in \Rn^{m\times n}$ contains independent normal distributed entries. In the experiment, we set $m=200, n=300, r=10$, $\gamma=10^{-3}$, $\alpha=\norm{{\X}^\star}_{\star}$, and $\beta=\frac{1}{2} \sum_{i,j\in \cI^c} (\hat{X}_{ij}^\star)^2$. 

We analyze the performance of the proposed algorithm in terms of normalized error defined as  $$\frac{\sum_{i,j\in \cI}(X_{ij}-M_{ij})^2}{\sum_{i,j\in \cI}(M_{ij})^2},$$ where $\X$ is the output generated by the algorithm. As for both the algorithms, the nuclear norm constraint \eqref{CV1} is always satisfied. We analyze the constraint violation for the constraint on unobserved entries given in \eqref{CV2}. The performance of the proposed algorithm is compared with \textbf{SCGD} \cite{mokhtari2018stochastic}. However, since \textbf{SCGD} cannot handle constraints in \eqref{CV2}, we incorporate the penalty on the unobserved entries as a regularization term. We remark that the matrix completion problem presented in \eqref{matrix_completion}-\eqref{CV2} is deterministic in nature, and the evaluation of the gradient $\nabla f(\X)$ requires access to all the observed entries in $\mathcal{I}$, which can be costly to compute at each iteration $t$. In such a scenario, a subset of the set $\mathcal{I}$ is often used as an unbiased estimate of the gradient \cite{hazan2016variance,mokhtari2018stochastic}. For instance, to obtain a stochastic approximation of gradient, mini batch method uses $b$ elements of $\mathcal{I}$ while the growing mini-batch
method uses a batch size of $b = \mathcal{O}(t^2)$ at step $t$ \cite{hazan2016variance}. In our experiments, to approximate the gradient at each instant $t$,  we use the average of stochastic gradients over time as given in \eqref{nh}  for \textbf{FW-CSOA}, and for \textbf{SCGD}, we use the average as suggested in \cite[Eq. (3)]{mokhtari2018stochastic}.

We run both algorithms for $T=3000$ iterations. For \textbf{FW-CSOA}, as dictated in theory, we set the parameters $\eta = \eta_0/T^{3/4}$, $\rho = \rho_0/\sqrt{T}$, and $\upsilon = \upsilon_0/\sqrt{T}$. Next, we tune the parameters $\eta_0, \rho_0$ and $\delta$ to achieve the best performance in terms of the  normalized error, resulting in $\eta_0 = 0.68, \rho_0=1.25$ and $\delta = 0.25$. Then, we tune $\upsilon_0$ to minimize the constraint violation without significantly sacrificing the optimality gap and finally set $\upsilon_0=0.77$. For \textbf{SCGD}, we choose a diminishing step size strategy as mentioned in \cite{mokhtari2018stochastic} and  set $\gamma_t=2/(t+8)$, $\rho_t=4/(t+8)^{2/3}$ while tuning the regularization parameter $\tau$ to obtain the best possible performance. Finally, we selected $\tau=5\times 10^{-6}$.

The evolution of the normalized error{,} as well as the constraint violation{,} is shown in Fig.~\ref{fig:matrix_comp}. For this experiment, we set the number of observed samples $b = 200$ for both algorithms. Fig.~\ref{NET} and Fig.\ref{CVT} illustrate the convergence path of normalized error and constraint violation, respectively. It can be observed that the proposed algorithm outperforms the compared method and is able to maintain almost zero average constraint violation as claimed. In the next experiment, we evaluate both algorithms for different values of $b$ and 
report the performance with respect to the total number of processed samples ($T\times b$) in Fig.~\ref{NES}. It can be observed that the performance of both the algorithms improves with the increase in {the} number of observed samples used at each iteration, which is the expected behavior. However, the proposed algorithm is able to achieve a small normalized error (Fig.~\ref{NES})  with almost negligible constraint violation{,} as shown in Fig.~\ref{CVS}.

\section{Conclusion} \label{sec:conclusion}
In this paper, we have considered a convex optimization problem where both the objective and constraint functions are in the form of expectations. We considered an augmented Lagrangian
relaxation of the problem and used a stochastic variant
of the primal-dual method to solve it. The proposed algorithm (\textbf{CSOA}) forces the constraint violation to be zero while maintaining the optimal convergence rate of $\mathcal{O}(T^{-1/2})$. We propose another algorithm \textbf{FW-CSOA} which is a projection-free version of \textbf{CSOA} and utilizes a momentum-based gradient tracking approach, achieving a rate of convergence of $\mathcal{O}(T^{-1/4})$ while maintaining zero average constraint violation. The proposed algorithms are tested on the fair classification and structural matrix completion problems, respectively, and were found to outperform state-of-the-art algorithms. As a future scope, it would be interesting to consider more general functional constraints and solve the problem using the conditional gradient method as discussed in \cite{lan2020conditional2}.

\appendices

\section{Proof of Corollary \ref{coro1}} \label{proof_corollary1}
Let $(\xu^\star,\lus)$ be the primal-dual optimal solution to \eqref{conservative_prob} while recall that $(\x^\star,\lam^\star)$ is the primal dual solution to \eqref{prob_form}. As $\upsilon \leq \frac{\sigma}{2}\leq \sigma$, there exists a strictly feasible primal vector $\tx$ such that $H_i(\tx)+\upsilon\leq H_i(\tx)+\sigma$ for $1\leq i \leq N$, and hence strong duality holds for \eqref{conservative_prob}. Therefore, we have:
\begin{align}
F(\xu^\star)&=\min_{\x} F(\x)+\ip{\lus,\H(\x)+\upsilon\one }\nonumber\\
& \leq  F(\x^\star)+\ip{\lus,\H(\x^\star)+\upsilon\one }\label{new_lag0}\\
& \leq F(\x^\star)+\upsilon\ip{\lus,\one } \label{new_lag}
\end{align}
where \eqref{new_lag0} follows from the optimality of $\xu^\star$ and \eqref{new_lag} from the fact that $\H(\x^\star)\leq 0$. Since $\tx$ is strictly feasible, we have from \eqref{new_lag0} that
\begin{align}
F(\xu^\star)&\leq F(\tx)+\ip{\lus,\H(\tx)+\upsilon\one}\nonumber\\& \leq F(\tx)+ (\upsilon-\sigma)\ip{\lus,\one}.
\label{new_lag2}
\end{align}
Equivalently we can write
\begin{align}
\one ^T\lus\leq \frac{F(\tx)-F(\xu^\star)}{\sigma-\upsilon} \leq \frac{2G_f \norm{\tx-\xu^\star}}{\sigma} \leq \frac{2G_fD}{\sigma}.
\label{new_lag3}
\end{align}
Here, in the second inequality we have used the fact that $\upsilon\leq\frac{\sigma}{2}$ and the last inequality comes from Assumption \ref{compact}. Now using \eqref{new_lag} and \eqref{new_lag3}, and using the fact that $ F(\x^\star) \leq F(\xu^\star)\leq F(\x^\star)+\upsilon\ip{\lus,\one }$, we get required expression
\begin{align}
|F(\xu^\star)-F(\x^\star)| \leq \frac{2G_fD}{\sigma} \upsilon := C\upsilon.
\end{align}

\section{Proof of Lemma \ref{lemma3}}\label{lemma3_proof}
The key to establishing Lemma \ref{lemma3} is Assumption \ref{smooth} regarding the smoothness of $f$ and $h_i$. Indeed, Assumption \ref{smooth} implies that $\L(\x,\lam,\th)$ is also smooth in $\x$ with the smoothness parameter depending on $\lam$. This can be seen by observing that
\begin{align}
&\|\nx\L(\x,\lam,\th)-\nabla_{\y}\L(\y,\lam,\th)\|
\nonumber\\&=\|\nabla f(\x,\th) \!-\!\nabla f(\y,\th)\!+\!\ip{\lam,\left(\nabla\h(\x,\th)-\nabla\h(\y,\th)\right)}\|\nonumber\\
& \leq (L_f+L_h(\one^\top\lam))\norm{\x-\y} \nonumber\\
& \leq (L_f+L_h\norm{\lam}_1)\norm{\x-\y} \label{laglip}
\end{align}
where we have used the triangle inequality. In other words, $\L(\x,\lam,\th)$ is $L_\lam$-smooth in $\x$ for $L_\lam = L_f + L_h\norm{\lam}_1$. 

Having established the smoothness condition, the rest of the proof follows along similar lines as that of Lemma 2 in \cite{cutkosky2019momentum}, with some modifications. Recall that $\L(\x,\lam) = \mbE_{\th}\left[\L(\x,\lam,\th)\right]$. We begin with introducing $(1-\rho)\nx\L(\x_{t-1},\lam_{t-1})$ and $\nx \L(\x_t,\lam_t)$ into \eqref{nh} to obtain
\begin{align}
\d_t - \nx\L(\x_t,\lam_t)&= (1-\rho)(\d_{t-1}-\nx\L(\x_{t-1},\lam_{t-1})) -(1-\rho)(\nx\L(\x_{t-1},\lam_{t-1},\th_t)-\nx\L(\x_{t-1},\lam_{t-1})) \nonumber\\&\quad+\nx\L(\x_t,\lam_t,\th_t)-\nx\L(\x_t,\lam_t). \label{dtupdate}
\end{align}
Let us consider the random variables $\sX_u = \nx \L(\x_u,\lam_u,\th_t)$ for $u \leq t$. Since $\th_t$ is independent of $\{\lam_u,\x_u\}_{u=1}^t$, it holds that
$\mbE\sX_u = \mbE[\nx \L(\x_u,\lam_u,\th_t)\mid\F_t] = \nx \L(\x_u,\lam_u)$ for all $u \leq t$. Taking full expectation, we therefore have that $\mbE\sX_u = \mbE[\nx \L(\x_u,\lam_u)]$ for all $u \leq t$. Taking squared norm on both sides of \eqref{dtupdate} and taking expectation, we thus have that
\begin{align}
\mbE\|\d_t-&\nx\L(\x_t,\lam_t)\|^2\nonumber\\&=(1-\rho)^2\mbE \|\d_{t-1}-\nx\L(\x_{t-1},\lam_{t-1})\|^2 + \mbE\|(1-\rho)(\sX_{t-1}-\mbE\sX_{t-1})+\sX_t-\mbE\sX_t\|^2.
\label{new_eq}
\end{align}
since the term $\d_{t-1}-\nx\L(\x_{t-1},\lam_{t-1})$ is independent of $\th_t$. The second term on the right of \eqref{new_eq} can be bounded as
\begin{align}\label{x1}
&\mbE\norm{(1-\rho)(\sX_{t-1}-\mbE\sX_{t-1})+\sX_t-\mbE\sX_t}^2
\nonumber\\& =\mbE\norm{(1-\rho)(\sX_{t-1}-\mbE\sX_{t-1}-\sX_t + \mbE\sX_t)+\rho(\sX_t-\mbE\sX_t)}^2\nonumber\\
&\leq 2(1-\rho)^2\mbE\norm{\sX_{t-1}-\sX_t}^2 + 2\rho^2\mbE\norm{\sX_t-\mbE\sX_t}^2
\end{align}
where we have used a norm inequality and the non-negativity of the variance. We can bound the variance of $\sX_t$ from Assumption \eqref{grad_norm_sq_x_zero} so that $\mbE\norm{\sX_t-\mbE\sX_t}^2 \leq 4B^2(1+\mbE\norm{\lam_t}^2)$. Now to bound first term of RHS of \eqref{x1} we introduce $\nx\L(\x_{t-1},\lam_{t},\th_t)$ in side the norm square, that is:
\begin{align}\label{x2}
&\mbE\norm{\sX_t-\sX_{t-1}}^2 \nonumber\\&=\mbE\norm{\sX_t\!-\sX_{t-1}\!+\nx\L(\x_{t-1},\lam_{t},\th_t)\!-\nx\L(\x_{t-1},\lam_{t},\th_t)}^2 \nonumber\\
&\leq 2\mbE\norm{\sX_t-\nx\L(\x_{t-1},\lam_{t},\th_t)}^2 +2\mbE\norm{\nx\L(\x_{t-1},\lam_{t},\th_t)-\sX_{t-1}}^2.
\end{align}
Recall from \eqref{laglip} that 
\begin{align}
2\mbE\|\sX_{t}-\nx\L(\x_{t-1},\lam_{t},\th_t)\|^2&\leq 2\norm{\x_t-\x_{t-1}}^2\mbE[L_{\lam_{t}}^2]\nonumber\\
&=2\eta^2\norm{\s_t-\x_t}^2\mbE[L_{\lam_{t}}^2]\leq2 \eta^2D^2\mbE[L_{\lam_{t}}^2]
\end{align}
where we have used the form of the primal update and Assumption \ref{compact}. Further, we can bound the second term of \eqref{x2} as:
\begin{align}\label{using lam update}
2\mbE\norm{\nx\L(\x_{t-1},\lam_t,\th_t)-\sX_{t-1}}^2& 
= 2\mbE\norm{\ip{\lam_t-\lam_{t-1},\nx \h(\x_{t-1},\th_t)}}^2\nonumber\\& 
\leq  2\mbE[\norm{\lam_t-\lam_{t-1}}^2\norm{\nx \h(\x_{t-1},\th_t)}^2]\nonumber\\& 
\leq  2\mbE\norm{\lam_t-\lam_{t-1}}^2N\sigma_{h}^2\nonumber\\& 
\leq  2N\sigma_{h}^2\eta^2\mbE\norm{\h(\x_{t-1},\th_t)+\textbf{1}\upsilon}^2\nonumber\\& 
\leq  4N\sigma_{h}^2\eta^2(\mbE\norm{\h(\x_{t-1},\th_t)}^2+N\upsilon^2)\nonumber\\& 
\leq  4N\sigma_{h}^2\eta^2(N\sigma_{\lam}^2+N\upsilon^2) 
\leq  8\eta^2N^2\sigma_{h}^2\sigma_{\lam}^2.
\end{align}
Here, in second inequality we used Cauchy Schwarz inequality, third inequality is obtained using  assumption \ref{boundedgrad}. In the fourth inequality, we have used the dual update equation of Algorithm 2 and next inequality is obtained from assumption \ref{boundedh}. Finally in last inequality we have used the fact that $\upsilon<\sigma_\lam$.
Now using bounds from \eqref{x1}-\eqref{using lam update} in \eqref{new_eq} we conclude Lemma \ref{lemma3} as:
	\begin{align}
\mbE\norm{\d_t - \nx\L(\x_{t},\lam_{t})}^2&\leq 8\rho^2 B^2(1+\mbE\norm{\lam_t}^2)+(1-\rho)^2\mbE \|\d_{t-1}-\nx\L(\x_{t-1},\lam_{t-1})\|^2\nonumber\\
&\quad +2(1-\rho)^2(\eta^2D^2\mbE[L_{\lam_t}^2]+8\eta^2N^2\sigma_{h}^2\sigma_\lambda^2)
\end{align}
\section{Proof of Lemma \ref{lemma4}} \label{lemma4_proof}
As in the proof of Lemma \ref{lemma1}, we divide this proof into two parts, establishing an upper bound on 
$\sum_{t=1}^T\mbE[\L (\x_t,\lam_t) -\L(\x, \lam_t)] $ and a lower bound on $\sum_{t=1}^T \mbE[\L(\x_t,\lam_t) -\L(\x_t, \lam)]$. We will use lower bound on $\sum_{t=1}^T \mbE[\L(\x_t,\lam_t) -\L(\x_t, \lam)] $ directly from \eqref{second half_lemma1} as the dual update is same in both the algorithms. Taking total expectation of \eqref{second half_lemma1} and summing it from  $t=1,2,\cdots,T$, we obtain
	%
	\begin{align}\label{second half_lemma1_2}
		\sum_{t=1}^T\mathbb{E}\left[\L(\x_t,\lam_t) -\L(\x_t, \lam)\right]
		\geq \frac{1}{2\eta}\mathbb{E}\left[\norm{\lam_{T+1}-\lam}^2-\norm{\lam_1 - \lam}^2 \right]  -\frac{\eta}{2} \sum_{t=1}^T\mathbb{E}\|\nabla_{\lam}\L(\x_t,\lam_t,\th_t)\|^2.
	\end{align} 
	By replacing the positive terms on the RHS of \eqref{second half_lemma1_2} by its lower bound $0$, we could write
	\begin{align}\label{second half_lemma1_3}
		\sum_{t=1}^T\mathbb{E}\left[\L(\x_t,\lam_t) -\L(\x_t, \lam)\right]
		\geq -\frac{1}{2\eta}\mathbb{E}[\norm{\lam_1 - \lam}^2] -\frac{\eta}{2} \sum_{t=1}^T\mathbb{E}\|\nabla_{\lam}\L(\x_t,\lam_t,\th_t)\|^2.
	\end{align} 
	Utilizing the upper bound mentioned in Corollary \ref{cor-bounds} (cf. \eqref{grad_norm_sq_lam_zero}) into the right hand side of \eqref{second half_lemma1_3}, we obtain
	\begin{align}\label{2second half_lemma1}
		\sum_{t=1}^T \mbE[\L(\x_t,\lam_t) -\L(\x_t, \lam)]\geq \!-\frac{1}{2\eta}\mbE\norm{\lam_{1}\!-\!\lam}^2   \!-\!\frac{\eta}{2}\!\sum_{t=1}^T (4N\sigma_{\lam}^2\!+\!4N\upsilon^2\!+\!2\delta^2\eta^2 \mbE\norm{\lam_t}^2)
	\end{align}
	By selection the parameter $\upsilon$ such that $\upsilon < \sigma_\lam$, we could further simplify the inequality in \eqref{2second half_lemma1} as 
	\begin{align}\label{2second half_lemma1_1}
		\sum_{t=1}^T \mbE[\L(\x_t,\lam_t) -\L(\x_t, \lam)]\geq -\frac{1}{2\eta}\mbE\norm{\lam_{1}-\lam}^2  -\frac{\eta}{2}\sum_{t=1}^T (8N\sigma_{\lam}^2+2\delta^2\eta^2 \mbE\norm{\lam_t}^2).
\end{align}
Next, we find an upper bound on $\sum_{t=1}^T\mbE[\L(\x_t,\lam_t)-\L(\x,\lam_t)]$. Since $\L(\x,\lam_t)$ is $L_{\lam_t}$-smooth, the quadratic upper bound implies that
\begin{align}\label{2ini_bound0}
	\L(\x_{t+1},\lam_t)& \leq  \L(\x_t,\lam_t)+\ip{\nx \L(\x_t,\lam_t), \x_{t+1}-\x_t }+ \frac{L_{\lam_t}}{2}\norm{\x_{t+1}-\x_t}^2\nonumber\\
	& = \L(\x_t,\lam_t)+\eta \ip{\nx \L(\x_t,\lam_t), \s_t-\x_t } + \frac{\eta^2 L_{\lam_t}}{2}\norm{\s_t-\x_t}^2 
\end{align}
where we have used the update equation \eqref{iterate_primal_FW}. Note that since $\s_t =\arg\min_{\s} \ip{\s,\d_t}$, we have that $\ip{\d_t,\s_t} \leq \ip{\d_t,\x}$ for any $\x\in\cX$. Therefore adding the non-negative term $\ip{\d_t,\x-\s_t}$ as well as adding and subtracting  $\eta\ip{\nx \L(\x_t,\lam_t),\x}$ to the right, we obtain
\begin{align}
	&\L(\x_{t+1},\lam_t) \leq  \L(\x_t,\lam_t) + \eta \ip{\nx \L(\x_t,\lam_t) - \d_t, \s_t-\x}+\eta \ip{\nx \L(\x_t,\lam_t), \x-\x_t} + \frac{\eta^2 L_{\lam_t}}{2}\norm{\s_t-\x_t}^2 \label{fwproofquad}
\end{align}

Noting that $\cX$ is compact from Assumption \ref{compact} and using the convexity of $\L(\x,\lam_t)$ with respect to $\x$, we obtain 
\begin{align}\label{2ini_bound}
	\L(\x_{t+1},\lam_t)&-\L(\x_t,\lam_t) \leq \eta D\|\nx \L(\x_t,\lam_t)-\d_t\| +\eta  (\L(\x,\lam_t)-\L(\x_t,\lam_t))+ \frac{\eta^2L_{\lam_t}D^2}{2}.
\end{align}
The left-hand side of \eqref{2ini_bound} can be lower bounded as:
\begin{align}\label{bound0}
	\L(\x_{t+1},\lam_t)-\L(\x_t,\lam_t)
	&=F(\x_{t+1})-F(\x_t)-\ip{\lam_t,\H(\x_t)-\H(\x_{t+1})}\nonumber\\
	&\geq F(\x_{t+1})-F(\x_t)-\frac{\eps_2\norm{\lam_t}^2}{2}-\frac{\norm{\H(\x_t)-\H(\x_{t+1})}^2}{2\eps_2}\nonumber\\
	&\geq F(\x_{t+1})-F(\x_t)-\frac{\eps_2\norm{\lam_t}^2}{2}-\frac{NG^2_h\eta^2D^2}{2\eps_2}
\end{align}
where we have used the Peter-Paul inequality for some $\eps_2>0$ which will be specified later, the Lipschitz continuity of $\{H_i\}_{i=1}^{N}$, and the compactness of $\cX$ which implies that
$\norm{\x_{t+1}-\x_t}^2=\eta^2\norm{\s_t-\x_t}^2\leq \eta^2D^2$. Substituting \eqref{bound0} into \eqref{2ini_bound}, taking expectation, and rearranging, we obtain
\begin{align}\label{exp1}
	\eta\mbE[\L(\x_t,\lam_t)-\L(\x,\lam_t)]&\leq \mbE[F(\x_t)-F(\x_{t+1})]+\frac{\eps_2\mbE\norm{\lam_t}^2}{2}+\frac{\eta^2NG^2_hD^2}{2\eps_2}\nonumber\\&\quad +\eta D \mbE\|\nx \L(\x_t,\lam_t)-\d_t\| +\frac{\eta^2 \mbE [L_{\lam_t}]D^2}{2}\nonumber\\&\leq \mbE[F(\x_t)-F(\x_{t+1})]+\frac{\eps_2\mbE\norm{\lam_t}^2}{2}+\frac{\eta^2NG^2_hD^2}{2\eps_2}\nonumber\\&\quad +\frac{\eps_1}{2}\mbE\norm{\nx \L(\x_t,\lam_t) - \d_t}^2 + \frac{\eta^2D^2}{2\eps_1} +\frac{\eta^2 \mbE [L_{\lam_t}]D^2}{2}
\end{align}
where in last expression we have used Peter-Paul inequality for some $\eps_1 > 0$ which will be specified later as well as the compactness of the set $\mathcal{X}$.

To further simplify we define $L:=\max\{L_f,L_h\sqrt{N},1\}$ and observe that
\begin{align}\label{here}
	\mbE[L_{\lam_t}]&=L_f + L_h \mbE\norm{\lam_t}_1\leq L_f + L_h \sqrt{N}\mbE\norm{\lam_t}\nonumber\\&\leq L(1+\mbE\norm{\lam_t}). 
\end{align} 
The bound in \eqref{here} also implies that 
\begin{align}\label{here2}
	\mbE[L_{\lam_t}] \leq \mbE[L_{\lam_t}^2] =L^2(1+\mbE\norm{\lam_t})^2
\end{align} 
Further, we can bound term $(1+\mbE\norm{\lam_t})^2$ as:
\begin{align}\label{here3}
(1+\mbE\norm{\lam_t})^2\leq 3+\frac{3}{2}\mbE\norm{\lam_t}^2,
\end{align}
 where we have used the Peter-Paul inequality and the fact that $(\mbE\norm{\lam_t})^2\leq \mbE\norm{\lam_t}^2$. 
 
  Now defining $\psi_t=\mbE[\L(\x_t,\lam_t)-\L(\x,\lam_t)]$,  $\phi_t=\norm{\d_t-\nx \L(\x_t,\lam_t)}^2$, and using \eqref{here2} and \eqref{here3} into \eqref{exp1} we can write:
\begin{align}\label{psi}
	\eta\psi_t&\leq \mbE[F(\x_t)-F(\x_{t+1})]+\frac{\eps_2\mbE\norm{\lam_t}^2}{2}+\frac{\eta^2NG^2_hD^2}{2\eps_2}\nonumber\\&\quad +\frac{\eps_1}{2}\phi_t + \frac{\eta^2D^2}{2\eps_1} +\frac{3\eta^2D^2L^2}{2}\left(1+\frac{\mbE\norm{\lam_t}^2}{2}\right).
\end{align}
Next, expressing Lemma \ref{lemma3} for $t=t+1$, summing the resulting expression over $t$, and using bound $1-\rho\leq 1$, we obtain:
\begin{align}\label{phit+1}
	\rho\sum_{t=1}^T\phi_{t}&\leq\phi_{1}+16\eta^2N^2T\sigma_{h}^2\sigma_{\lam}^2+8\rho^2 B^2T+6\eta^2D^2L^2T\nonumber\\
	&\quad+8\rho^2 B^2\sum_{t=1}^T\mbE\norm{\lam_{t+1}}^2
	+3\eta^2D^2L^2\sum_{t=1}^T\mbE\norm{\lam_{t+1}}^2.
\end{align}
Bounding the first term of RHS of the dual update equation (step 6) of Algorithm \ref{algo_2} as $(1-\eta^2\delta) \leq 1$ and taking norm square of the resulting expression, we get: 
\begin{align}\label{lam}
	\mbE\norm{\lam_{t+1}}^2&
	\leq \mbE\norm{\lam_t+\eta(\h(\x_t,\th_t)+\textbf{1}\upsilon)}^2\nonumber\\&\leq 2\mbE\|\lam_t\|^2+2\eta^2\mbE\|\h(\x_t,\th_t)+\textbf{1}\upsilon)\|^2\nonumber\\&\leq 2\mbE\|\lam_t\|^2+4\eta^2\mbE\|\h(\x_t,\th_t)\|^2+4\eta^2N\upsilon^2\nonumber\\&\leq 2\mbE\|\lam_t\|^2+4\eta^2N\sigma_{\lam}^2+4\eta^2N\upsilon^2\nonumber\\&\leq 2\mbE\|\lam_t\|^2+8\eta^2N\sigma_{\lam}^2
\end{align}
here in second and third inequality we have used norm inequality while the fourth inequality is obtained using  Assumption \ref{boundedh} and the last inequality is obtained using the fact that $\upsilon<\sigma_\lam$. Now substituting \eqref{lam} into \eqref{phit+1}, we obtain
\begin{align}\label{phi}
	\rho\sum_{t=1}^T\phi_{t}&\leq\phi_{1}+16\eta^2N^2T\sigma_{h}^2\sigma_{\lam}^2+8\rho^2 B^2T+6\eta^2D^2L^2T\nonumber\\&\quad+16\rho^2 B^2\sum_{t=1}^T\mbE\norm{\lam_{t}}^2
	+6\eta^2D^2L^2\sum_{t=1}^T\norm{\lam_{t}}^2+64\rho^2 \eta^2B^2N\sigma_{\lam}^2+24\eta^4D^2L^2N\sigma_{\lam}^2 
\end{align}
Summing \eqref{psi} over $t=1,2,\cdots,T$  and using \eqref{phi}, we get:
\begin{align}\label{psi+phi}
	\sum_{t=1}^T\psi_t&\leq \frac{1}{\eta}\mbE[F(\x_1)-F(\x_{T+1})]+\frac{\eps_2}{2\eta}\sum_{t=1}^T\mbE\norm{\lam_t}^2+\frac{\eta NG^2_hD^2T}{2\eps_2} + \frac{\eta D^2T}{2\eps_1}+\frac{3\eta D^2L^2T}{2} \nonumber\\&\quad+\frac{3\eta D^2L^2}{2}\sum_{t=1}^T\frac{\mbE\norm{\lam_t}^2}{2}+\frac{\eps_1 \phi_1}{2\rho\eta}+\frac{8\eta\eps_1T}{\rho}N^2\sigma_{h}^2\sigma_{\lam}^2+\frac{4\eps_1\rho B^2T}{\eta}+\frac{3\eps_1\eta}{\rho}D^2L^2T\nonumber\\&\quad+32\eps_1 \rho \eta B^2N\sigma_{\lam}^2+\frac{\eps_1\eta^3}{\rho}12D^2L^2N\sigma_{\lam}^2 +\frac{8\eps_1\rho B^2}{\eta}\sum_{t=1}^T\mbE\norm{\lam_{t}}^2+\frac{3\eps_1\eta}{\rho}D^2L^2 \sum_{t=1}^T\mbE\norm{\lam_{t}}^2.
\end{align}

Since $F$ is Lipschitz continuous, the first term on the right of \eqref{psi+phi} can be bounded by  $G_fD$. Further, $\phi_1=\mbE\|\d_{1}-\nx\L(\x_{1},\lam_1)\|^2\leq \sigma_f^2$ if we initialize $\d_{1}=\lam_{1}=\textbf{0}$. Using these bounds and subtracting \eqref{2second half_lemma1_1} from  \eqref{psi+phi} we obtain the desired result in  Lemma \ref{lemma4}, with $A_1$ and $A_2$ defined as:

\begin{align}\label{A_1}
	&A_1:=4\eta N\sigma_{\lam}^2T+\frac{G_fD}{\eta}+\frac{\eps_1\sigma_f^2}{2\eta\rho}+\frac{\eta NG^2_hD^2T}{2\eps_2}\\&\quad + \frac{\eta D^2T}{2\eps_1}+\frac{3\eta D^2L^2T}{2}+\frac{8\eta \eps_1T}{\rho}N^2\sigma_{h}^2\sigma_{\lam}^2+\frac{4\eps_1\rho B^2T}{\eta}\nonumber\\&\quad +  \frac{3\eps_1\eta}{\rho}D^2L^2T+32\eps_1 \rho \eta B^2N\sigma_{\lam}^2+\frac{\eps_1\eta^3}{\rho}12D^2L^2N\sigma_{\lam}^2 \nonumber \\&\label{A_2}
	A_2:=\frac{\eps_2}{2\eta}+\frac{3\eta D^2L^2}{4}+\frac{8\eps_1\rho B^2}{\eta}+\frac{3\eps_1\eta}{\rho}D^2L^2.
\end{align}
\section{Proof of Theorem \ref{theorem2}} \label{proof_throrem2}
Theorem \ref{theorem2} can be proved along the same line as Lemma \ref{lemma2} and Theorem \ref{theorem1}. We start with  taking expectation in \eqref{diff_eq} and summing it over $t=1,2,\cdots,T$. Then, using bound from Lemma \ref{lemma4} we can write:
\begin{align}\label{expected2}
	&\sum_{t=1}^T\!\mbE [F(\x_t)\!-\!F(\x)\! +\! \frac{\delta \eta}{2}\!\left( \norm{\lam_t}^2\!-\!\norm{\lam}^2\right)\!+\!\ip{\lam,\H(\x_t)\!+\!\upsilon\one }]\nonumber\\&\leq A_1+\frac{\norm{\lam}^2}{2\eta} +\frac{\eta}{2}\left(\frac{2A_2}{\eta}+2\delta^2\eta^2\right)\sum_{t=1}^T \mbE\norm{\bblam_t}^2
\end{align}
Rearranging, we obtain
\begin{align}\label{expected3}
	\sum_{t=1}^T\mbE[F(\x_t)-F(\x)]& - \left(\frac{\delta \eta T}{2}+\frac{1}{2\eta}\right) \norm{\lam}^2+\sum_{t=1}^T\mbE\left[\ip{\lam,\H(\x_t)+\upsilon\one }\right]\nonumber\\&\leq A_1+ \frac{\eta}{2} \left(\frac{2A_2}{\eta}+2\delta^2\eta^2-\delta \right)\sum_{t=1}^T \mbE \norm{\bblam_t}^2
\end{align}
Let us concentrate on the term $\tfrac{2A_2}{\eta}+2\delta^2\eta^2-\delta$ which is quadratic in $\delta$. For roots of this quadratic equation to be real we must have $1-4(2\eta^2)(\frac{2A_2}{\eta})=1-16\eta A_2\geq 0$, that is:  
\begin{align}\label{condtn}
1\!-\!16\left(\frac{\eps_2}{2}\!+\!\frac{3\eta^2 D^2L^2}{4}\!+\!8\eps_1\rho B^2\!+\!\frac{3\eps_1\eta^2}{2\rho}D^2L^2\right)\!\geq\!0
\end{align}
where, we have substituted value of $A_2$ from \eqref{A_2}.
Now rearranging \eqref{condtn} we obtain:
\begin{align}
\eta^2\leq \frac{1-16\left(\frac{\eps_2}{2}+8\eps_1\rho B^2\right)}{16\left(\frac{3 L^2 D^2}{4}+\frac{3\eps_1L^2 D^2 }{\rho}\right)}
\end{align}
Substituting $\eps_1=\frac{T^{-a}}{16B}$, $\rho=\frac{T^{-a}}{8B}$ and $\eps_2=\frac{T^{-b}}{8}$, where $a$ and $b$ are positive scalar we get,  \begin{align} \label{z1}
\eta^2\leq \frac{1-T^{(-2a-b)}}{36L^2D^2}
\end{align}
 For $\eta$ defined in \eqref{z1}, we can set  $\delta=\frac{9L^2D^2}{1-T^{(-2a-b)}}$ so that $\frac{2A_2}{\eta}+2\delta^2\eta^2\leq \delta$. For simplicity, we write the denominator of $\delta$ for large $T$ as $1-T^{(-2a-b)}\approx 1$ and hence set  $\delta=9L^2D^2$.
This allows us to drop the terms containing $\mbE\norm{\lam_t}^2$ on the right of \eqref{expected3}. Therefore, we have: 
\begin{align}\label{expected4}
	&\sum_{t=1}^T\mbE[F(\x_t)-F(\x)] - \left(\frac{\delta \eta T}{2}+\frac{1}{2\eta}\right) \norm{\lam}^2+\sum_{t=1}^T\mbE\left[\ip{\lam,\H(\x_t)+\upsilon\one }\right]\leq A_1\end{align}
It is then possible to tighten the bound in \eqref{expected4} in the same way as we did in \eqref{opt_grad}. Proceeding in a similar fashion as in \eqref{main_equation_expctd8}-\eqref{use_coro},  we get following three bounds:
\begin{align}
	&\!\sum_{t=1}^T \mbE[F(\x_t)\!-\!F(\x)]\!+\!		\frac{\left[\sum_{t=1}^T\mbE[\H(\x_t)+\upsilon\one ]\right]_{+}^2}{2(\delta \eta T +\frac{1}{\eta})}\leq  A_1 \label{main_equation_expctd82}\\
	&\sum_{t=1}^T\mbE[F(\x_t)-F(\x_{\upsilon}^\star)]\leq A_1
	\label{bound_on_err_seq2}\\
	&\sum_{t=1}^T\mbE[F(\x_t)-F(\x^\star)]\leq A_1+ C \upsilon  T.
	\label{bound_on_err_seq222}
\end{align}
For constraint violation, we will consider standard Lagrangian the same we did in \eqref{main_equation_expctd443}  and bound difference $\mbE[ \cL(\x_t,\lam) -\cL(\x, \lam_t)] $ using bound from Lemma \ref{lemma4} as follows:
\begin{align}\label{th2:21}
	&\sum_{t=1}^T\mbE[ \cL(\x_t,\lam) -\cL(\x, \lam_t)]\nonumber\\&\leq A_1+\frac{1}{2\eta} \norm{\lam}^2+\sum_{t=1}^T \frac{\eta}{2}\left(\frac{2A_2}{\eta}+2\delta^2\eta^2\right) \mbE \norm{\bblam_t}^2+\sum_{t=1}^T\frac{\delta \eta}{2}\norm{\lam}^2-\sum_{t=1}^T\frac{\delta \eta}{2}\mbE\norm{\lam_t}^2\nonumber\\&=  A_1\!+\! \frac{1}{2\eta}\mbE\norm{\lam}^2\!+\!\sum_{t=1}^T\frac{\eta}{2} \left(\frac{2A_2}{\eta}+2\delta^2\eta^2-\delta \right) \mbE \norm{\bblam_t}^2+\frac{\delta \eta T}{2}\norm{\lam}^2
\end{align}
Setting $\delta=9L^2D^2$, we can  drop the terms multiplying $\mbE\norm{\lam_t}^2$ as we did in \eqref{expected3}. Thus, we have:
\begin{align}\label{th2:22}
	\sum_{t=1}^T&\mbE[ \cL(\x_t,\lam)\! -\!\cL(\x, \lam_t)]\leq A_1\! +\! \left(\frac{1}{2\eta}\!+\!\frac{\delta \eta T}{2} \right)\norm{\lam}^2
\end{align}
Now using \eqref{stochastic_lagrangian2233} and using the bound from \eqref{th2:22} with $\x = \xus$ and $\lam = \lus + \one_i$, we obtain
\begin{align} \label{th22}
	\sum_{t=1}^T \E{H_i(\x_t)}+\upsilon T &\leq  \sum_{t=1}^T\E{\cL(\x_t,\one _i+\lus)-\cL(\xus,\lam_t)}\nonumber\\
	&\leq A_1 + \left({1}/{2\eta}+{\delta \eta T}/{2}\right)\norm{\lus+\one_i}^2\nonumber\\
	&\leq A_1 + \left({1}/{2\eta}+{\delta \eta T}/{2}\right)\left( 1 + C^2\right)
\end{align}
Therefore, in order to ensure that the constraint violation is less than or equal to zero, we choose 
\begin{align}\label{ups}
	\upsilon = \frac{A_1}T + \left(\frac{1}{2\eta T}+\frac{\delta \eta }{2} \right)\left( 1 + C^2\right)
\end{align}
which implies that the optimality gap is bounded by 
\begin{align}\label{gap2}
	\sum_{t=1}^T& \E{F(\x_t)}-F(\x^\star)   \leq  A_1+ C\left(A_1 + \left(\frac{1}{2\eta}+\frac{\delta \eta T}{2} \right)\left( 1 + C^2\right)\right)
\end{align}

Now setting $\eps_1=\rho=\frac{T^{-a}}{8B}$ and $\eps_2=\frac{T^{-b}}{8}$ in \eqref{A_1}, the bound in \eqref{gap2} can be minimized by setting $a=b=0.5$ and $\eta=\tfrac{T^{-3/4}}{6LD}=\tfrac{\hat{C_1}}{T^{3/4}}$. Substituting these values in \eqref{ups} and \eqref{gap2}, we obtain:
\begin{subequations}\label{final_re2}
	\begin{align}&A_1\leq A T^{3/4},\\
		&\upsilon=\hat{C_2}T^{-1/4} \;\text{and} \\
		&\frac{1}T\sum_{t=1}^T \E{F(\x_t)}\!-\!F(\x^\star) \leq  \frac{\hat{K}}T(T^{3/4})\!=\!\hat{K}(T^{-1/4}),
	\end{align}
\end{subequations}
where
\begin{subequations}\label{final_cons}
	\begin{align}
		&A=\tfrac{LD}{16}\left(96G_f+24\sigma_{f}^2+11\right)+\tfrac{N\sigma_{\lam}^2}{3LD}\left(\tfrac{17}{24}+2N\sigma_{h}^2\right)\nonumber\\&\quad\quad+\tfrac{1}{3L}\left(NG_h^2D+4DB\right) \\& \hat{C_2}=A+(\tfrac{15LD}{4})(1+C^2),\\& \hat{K}=A+C(A+(\tfrac{15LD}{4})(1+C^2))
	\end{align}
\end{subequations}
which implies the required result. 

\section{Proof of Proposition \ref{prop1} (Lower Bound on Optimality Gap):}
\label{bound_dual_optimal}

(a) From \eqref{new_lag3} and using the fact that $\norm{\y}\leq \norm{\y}_1=\textbf{1}^T\y$ for any positive vector $\y \geq 0$, we can write:
\begin{align}\label{bound_on_lam2}
\norm{\lus}^2\leq\left(\frac{2G_f D}{\sigma}\right)^2=C^2
\end{align}
Let us consider a set $\mathcal{S}$ as:
\begin{align}\label{optimal set}
\mathcal{S}:=\left\{\lam \geq 0 | \norm{\lam} \leq C+r\right\}
\end{align} 
where $C=\frac{2G_f D}{\sigma}$ with a scalar $r>0$. Note that the set $\mathcal{S}$ contains the set of
dual optimal solutions. For any $\lam \in \mathcal{S}$ from the dual update expression \eqref{iterate_dual},  we can write,
\begin{align}
\norm{\lam_{t+1}-\lam}^2=\norm{[\lam_t+\eta\nabla_{\lam}\L(\x_t,\lam_t,\th_t)]_{+}-\lam}^2.
\label{temp:lam1}
\end{align}
Using the non-expansiveness property of the projection operator $[\cdot]_{+}$ and expanding the square, we obtain:
\begin{align}
\norm{\lam_{t+1}-\lam}^2\leq &\norm{\lam_t+\eta\nabla_{\lam}\L(\x_t,\lam_t,\th_t)-\lam}^2\nonumber\\
=&\norm{\lam_t-\lam}^2+2\eta\ip{\nabla_{\lam}\L(\x_t,\lam_t,\th_t),\lam_t-\lam} +\eta^2\norm{\nabla_{\lam}\L(\x_t,\lam_t,\th_t)}^2.\label{above1}
\end{align} 
Reordering the terms in \eqref{above1} and taking total expectation we get
\begin{align}\label{temp:grad transpose_lam21}
\mbE[\ip{\lam-\lam_t, \nabla_{\lam}\L (\x_t,\lam_t)}]\leq \frac{1}{2\eta}\mbE\left[\norm{\lam_t-\lam}^2-\norm{\lam_{t+1}-\lam}^2\right]+\frac{\eta}{2}\mbE \norm{\nabla_{\lam}\L_t (\x_t,\lam_t)}^2.
\end{align}
Note that the Lagrangian is concave with respect to $\lam$, it holds that 
\begin{align}
\mbE[\ip{\lam_t-\lam^\star_{\upsilon},\nabla_{\lam}\L (\x_t,\lam_t)}]&\leq \mbE[\L(\x_t,\lam_t)] -\mbE[\L(\x_t, \lam^\star_{\upsilon})]\nonumber\\&=\mbE[\L(\x_t,\lam_t)] -\mbE[\mathcal{L}(\x_t, \lam^\star_{\upsilon})]+\frac{\delta \eta}{2} \norm{\lam^\star_{\upsilon}}^2\nonumber\\&\leq \mbE[\L(\x_t,\lam_t)] -\mbE[\mathcal{L}(\x^\star_{\upsilon}, \lam_t)]+\frac{\delta \eta}{2}\norm{\lam^\star_{\upsilon}}^2\nonumber\\&= \mbE[\L(\x_t,\lam_t)] -\mbE[\mathcal{L}(\x^\star_{\upsilon}, \lam_t)]+\frac{\delta \eta}{2}\norm{\lam^\star_{\upsilon}}^2+\frac{\delta \eta}{2}\mbE\norm{\lam_t}^2-\frac{\delta \eta}{2}\norm{\lam_t}^2\nonumber\\&= \mbE[\L(\x_t,\lam_t)] -\mbE[\L(\x^\star_{\upsilon}, \lam_t)]+\frac{\delta \eta}{2}\norm{\lam^\star_{\upsilon}}^2-\frac{\delta \eta}{2}\mbE\norm{\lam_t}^2
\label{temp:intro strng_lam1}
\end{align}
where, in second inequality we have expressed augmented lagrangian in terms of standard lagrangian, while in third inequality we have used relation from \eqref{lag_bound_trick} and in last inequality we express back in terms of augmented lagrangian. Now using \eqref{first half_lemma1} we can bound first two term of RHS of \eqref{temp:intro strng_lam1}. Thus we can write:
\begin{align}\label{new1}
\mbE[\ip{\lam_t-\lam^\star_{\upsilon},\nabla_{\lam}\L (\x_t,\lam_t)}]&\leq \frac{1}{2\eta} \mbE \left[ \norm{\x_t-\x^\star_{\upsilon}}^2  -  \norm{\x_{t+1} - \x^\star_{\upsilon}}^2\right]+\frac{\eta}{2}\mbE\norm{ \nabla_{\x}\L_t(\x_t,\lam_t)}^2+\frac{\delta \eta}{2}\norm{\lam^\star_{\upsilon}}^2-\frac{\delta \eta}{2}\mbE\norm{\lam_t}^2
\end{align}
Next, consider that 
\begin{align}\label{new2}
\mbE[\ip{\lam-\lam^\star_{\upsilon},\nabla_{\lam}\L (\x_t,\lam_t)}]&=\mbE[\ip{\lam-\lam^\star_{\upsilon}-\lam_t+\lam_t,\nabla_{\lam}\L (\x_t,\lam_t)}]\nonumber\\&=\mbE[\ip{\lam-\lam_t,\nabla_{\lam}\L(\x_t,\lam_t)}]+\mbE[\ip{\lam_t-\lam^\star_{\upsilon},\nabla_{\lam}\L(\x_t,\lam_t)}] 
\end{align}
Using bounds from \eqref{temp:grad transpose_lam21} and \eqref{new1} in \eqref{new2} we have:

\begin{align}\label{new3}
\mbE[\ip{\lam-\lam^\star_{\upsilon},\nabla_{\lam}\L (\x_t,\lam_t)}]&\leq \frac{1}{2\eta}\mbE\left[\norm{\lam_t-\lam}^2-\norm{\lam_{t+1}-\lam}^2\right]+\frac{\eta}{2}\mbE \norm{\nabla_{\lam}\L_t (\x_t,\lam_t)}^2\nonumber\\& \quad +\frac{1}{2\eta} \mbE \left[ \norm{\x_t-\x^\star_{\upsilon}}^2  -  \norm{\x_{t+1} - \x^\star_{\upsilon}}^2\right]+\frac{\eta}{2}\mbE\norm{ \nabla_{\x}\L_t(\x_t,\lam_t)}^2+\frac{\delta \eta}{2}\norm{\lam^\star_{\upsilon}}^2-\frac{\delta \eta}{2}\mbE\norm{\lam_t}^2\nonumber\\&\leq \frac{1}{2\eta}\mbE\left[\norm{\lam_t-\lam}^2-\norm{\lam_{t+1}-\lam}^2\right]+\frac{1}{2\eta} \mbE \left[ \norm{\x_t-\x^\star_{\upsilon}}^2  -  \norm{\x_{t+1} - \x^\star_{\upsilon}}^2\right]\nonumber\\& \quad +\frac{\eta}{2} \left[2B^2(1+\mbE\norm{\lam_t}^2) + 2N\sigma^2_{\lam} + 2\delta^2\eta^2\mbE\norm{\lam_t}^2\right]  +\frac{\delta \eta}{2}\norm{\lam^\star_{\upsilon}}^2-\frac{\delta \eta}{2}\mbE\norm{\lam_t}^2.
\end{align}
where in last inequality we have used bound from \eqref{grad_norm_sq_x_zero} and \eqref{grad_norm_sq_lam_zero}. Now evaluating the gradient $\nabla_{\lam}\L (\x_t,\lam_t)$ we can simplify LHS of \eqref{new3} as:
\begin{align}\label{new32}
\mbE[\ip{\lam-\lam^\star_{\upsilon},\nabla_{\lam}\L (\x_t,\lam_t)}]&=\mbE[\ip{\lam-\lam^\star_{\upsilon},\H(\x_t)+\upsilon \textbf{1}}-\delta \eta \ip{\lam-\lam^\star_{\upsilon},\lam_t}]\nonumber\\&\geq \mbE[\ip{\lam-\lam^\star_{\upsilon},\H(\x_t)+\upsilon \textbf{1}}]-\delta \eta \left[\norm{\lam-\lam^\star_{\upsilon}}^2+\frac{\mbE \norm{\lam_t}^2}{4}\right],
\end{align}
where, in last inequality we have used Young's inequality. Now using bound \eqref{new32} in \eqref{new3}, rearranging and defining $P:=2B^2+2N\sigma^2_{\lam} $ and $Q=2B^2$ we can write:
\begin{align}\label{new31}
\mbE[\ip{\lam-\lam^\star_{\upsilon},\H(\x_t)+\upsilon \textbf{1}}]-\delta \eta \norm{\lam-\lam^\star_{\upsilon}}^2&\leq\frac{1}{2\eta}\mbE\left[\norm{\lam_t-\lam}^2-\norm{\lam_{t+1}-\lam}^2\right]+\frac{1}{2\eta} \mbE \left[ \norm{\x_t-\x^\star_{\upsilon}}^2  -  \norm{\x_{t+1} - \x^\star_{\upsilon}}^2\right]\nonumber\\& \quad + \frac{\eta}{2}[P+(Q +2\delta^2\eta^2-\frac{\delta}{2})\mbE\norm{\lam_t}^2]  +\frac{\delta \eta}{2}\norm{\lam^\star_{\upsilon}}^2.
\end{align}
Now choosing $\delta$ such that  $2Q +4\delta^2\eta^2 \leq \delta$ will make the terms containing $\mbE\norm{\lam_t}^2$ negative, allowing us to drop them from the right, yielding:
\begin{align}\label{new4}
\mbE[\ip{\lam-\lam^\star_{\upsilon},\H(\x_t)+\upsilon \textbf{1}}]&\leq \frac{1}{2\eta}\mbE\left[\norm{\lam_t-\lam}^2-\norm{\lam_{t+1}-\lam}^2\right]+\frac{1}{2\eta} \mbE \left[ \norm{\x_t-\x^\star_{\upsilon}}^2  -  \norm{\x_{t+1} - \x^\star_{\upsilon}}^2\right]\nonumber\\&\quad 
+\frac{\delta \eta}{2}\norm{\lam^\star_{\upsilon}}^2+\delta \eta \norm{\lam-\lam^\star_{\upsilon}}^2.
\end{align}

Summing \eqref{new4} for $t=1,2,\cdots,T$, using the fact that $\lam_{1}=0$ and using the boundness assumption of set $\mathcal{X}$, we obtain: 
\begin{align}\label{new5}
\sum_{t=1}^{T}\mbE[\ip{\lam-\lam^\star_{\upsilon},\H(\x_t)+\upsilon \textbf{1}}]&\leq\frac{1}{2\eta}\norm{\lam_1 -\lam}^2+\frac{1}{2\eta}\norm{\x_1-\x^\star_{\upsilon}}^2
+\frac{\delta \eta T}{2}\norm{\lam^\star_{\upsilon}}^2+\delta \eta \norm{\lam-\lam^\star_{\upsilon}}^2\nonumber\\& \leq\frac{1}{2\eta}\norm{\lam}^2+\frac{D^2}{2\eta}
+\frac{\delta \eta T}{2}\norm{\lam^\star_{\upsilon}}^2+2\delta \eta T(\norm{\lam}^2+\norm{\lam^\star_{\upsilon}}^2).
\end{align}

Let us denote $b=\H(\x_t)+\upsilon \textbf{1}$ and $b^+=[\H(\x_t)+\upsilon \textbf{1}]_+$ and for $r>0$ define a dual vector $\hat{\lam}$ as follows:

\begin{align}\label{new_lam}
\hat{\lam}=\lam^\star_{\upsilon}+r\frac{b^+}{\norm{b^+}}
\end{align}

Note since, $\lus\geq 0$, $b^+\geq 0$ and $r>0$, we have $\hat{\lam}\geq 0$. Using bound from \eqref{bound_on_lam2}, we can bound $||\hat{\lam}||$ as:

\begin{align}
||\hat{\lam}||\leq ||\lam^\star_{\upsilon}||+r\leq C+r,
\end{align}

this implies that $\hat{\lam}\in \mathcal{S}$. Now, from \eqref{new5} we can write:
\begin{align}\label{new555}
\max_{\lam \in \mathcal{S}}\sum_{t=1}^{T}\mbE[\ip{\lam-\lam^\star_{\upsilon},\H(\x_t)+\upsilon \textbf{1}}]&\leq \frac{1}{2\eta}\max_{\lam \in \mathcal{S}}\norm{\lam}^2+\frac{D^2}{2\eta}
+\frac{\delta \eta T}{2}\norm{\lam^\star_{\upsilon}}^2+2\delta \eta T\left( \max_{\lam \in \mathcal{S}} \norm{\lam}^2+\norm{\lam^\star_{\upsilon}}^2\right)\nonumber\\&= \left(\frac{1}{2\eta}+2\delta \eta T\right)\max_{\lam \in \mathcal{S}}\norm{\lam}^2+\frac{D^2}{2\eta}
+\left(\frac{\delta \eta T}{2}+2\delta \eta T\right)\norm{\lam^\star_{\upsilon}}^2.
\end{align}

As we have $\hat{\lam}\in \mathcal{S}$, from definition of $b$ and \eqref{new555} we obtain:

\begin{align}\label{new100}
\sum_{t=1}^{T}\mbE[\ip{\hat{\lam}-\lam^\star_{\upsilon},b}]&=\sum_{t=1}^{T}\mbE[\ip{\hat{\lam}-\lam^\star_{\upsilon},\H(\x_t)+\upsilon \textbf{1}}]\nonumber\\&\leq \max_{\lam \in \mathcal{S}}\sum_{t=1}^{T}\mbE[\ip{\lam-\lam^\star_{\upsilon},\H(\x_t)+\upsilon \textbf{1}}]\nonumber\\&\leq \left(\frac{1}{2\eta}+2\delta \eta T\right)\max_{\lam \in \mathcal{S}}\norm{\lam}^2+\frac{D^2}{2\eta}
+\left(\frac{\delta \eta T}{2}+2\delta \eta T\right)\norm{\lam^\star_{\upsilon}}^2.
\end{align}

Also, since $\hat{\lam}-\lam^\star_{\upsilon}=r\frac{b^+}{\norm{b^+}}$, we have $\ip{\hat{\lam}-\lam^\star_{\upsilon},b}=r||b^+||$. Thus, from the definition of $b$ we obtain
\begin{align}
\sum_{t=1}^{T}\mbE[\ip{\hat{\lam}-\lam^\star_{\upsilon},b}]=r\sum_{t=1}^{T}\mbE||[\H(\x_t)+\upsilon \textbf{1}]_+||
\end{align}

Substituting the preceding equality in \eqref{new100} and dividing both sides by $r$ we obtain

\begin{align}\label{divide_r}
\sum_{t=1}^{T}\mbE||[\H(\x_t)+\upsilon \textbf{1}]_+||\leq \left(\frac{1}{2r\eta}+\frac{2\delta \eta T}{r}\right)\max_{\lam \in \mathcal{S}}\norm{\lam}^2+\frac{D^2}{2r\eta}
+\left(\frac{\delta \eta T}{2r}+\frac{2\delta \eta T}{r}\right)\norm{\lam^\star_{\upsilon}}^2.
\end{align}

From definition of set $\mathcal{S}$ we have,
\begin{align}\label{l1}
\max_{\lam \in \mathcal{S}}\norm{\lam}^2\leq (C+r)^2,
\end{align}
Using bounds from \eqref{bound_on_lam2} and  \eqref{l1} and substituting $\eta=\frac{\hat{C}}{\sqrt{T}}$ in \eqref{divide_r}, we obtain:
\begin{align}\label{divide_r2}
\sum_{t=1}^{T}\mbE||[\H(\x_t)+\upsilon \textbf{1}]_+||&\leq \left(\frac{1}{2r\eta}+\frac{2\delta \eta T}{r}\right)(C+r)^2+\frac{D^2}{2r\eta}
+\left(\frac{\delta \eta T}{2r}+\frac{2\delta \eta T}{r}\right)C^2\nonumber\\&=\sqrt{T}\left[\frac{(\hat{C}^2+4\delta)(C+r)^2+D^2\hat{C}^2+5\delta C^2}{2\hat{C}r}\right]\nonumber\\&=\sqrt{T} Q,
\end{align}
where $Q=\left[\frac{(\hat{C}^2+4\delta)(C+r)^2+D^2\hat{C}^2+5\delta C^2}{2\hat{C}r}\right]$. Now for the the optimal pair ($\x^\star_{\upsilon}, \lam^\star_{\upsilon}$) of the standard Lagrangian \eqref{lagrangian}  we can write:
\begin{align}
\mbE[F(\x_t)]&=\mbE[F(\x_t)+\ip{\lam^\star_{\upsilon},\H(\x_t)+\upsilon \textbf{1}}-\ip{\lam^\star_{\upsilon},\H(\x_t)+\upsilon \textbf{1}}]\nonumber\\&=\mbE[\mathcal{L}(\x_t-\lam^\star_{\upsilon})]-\mbE[\ip{\lam^\star_{\upsilon},\H(\x_t)+\upsilon \textbf{1}}]
\end{align}
From \eqref{lag_bound_trick} we have:
\begin{align}\label{lag_bound_trick2}
\mbE[\cL(\x,\lam^\star_{\upsilon})]\geq\mbE[\cL(\x^\star_{\upsilon},\lam^\star_{\upsilon})]=\mbE[F(\x^\star_{\upsilon})].
\end{align}
Hence,
\begin{align}\label{eq_zz}
\mbE[F(\x_t)]\geq \mbE[F(\x^\star_{\upsilon})]- \mbE[\ip{\lam^\star_{\upsilon},\H(\x_t)+\upsilon \textbf{1}}]	
\end{align}
Summing for $t=1,2,\cdots,T$ we get:
\begin{align}\label{eq_zz_sum}
\sum_{t=1}^{T}\mbE[F(\x_t)-F(\x^\star_{\upsilon})]\geq-\sum_{t=1}^{T} \mbE[\ip{\lam^\star_{\upsilon},\H(\x_t)+\upsilon \textbf{1}}]	
\end{align}
Also, since $\lam^\star_{\upsilon}\geq 0$ and $\mbE[\H(\x_t)+\upsilon \textbf{1}]\leq\mbE[\H(\x_t)+\upsilon \textbf{1}]_+$, it follows that
\begin{align}\label{eq_zz2}
-\sum_{t=1}^{T}\mbE[\ip{\lam^\star_{\upsilon},\H(\x_t)+\upsilon \textbf{1}}]	\geq-\sum_{t=1}^{T}\mbE[\ip{\lam^\star_{\upsilon},[\H(\x_t)+\upsilon \textbf{1}}]_+]	\geq-\norm{\lam^\star_{\upsilon}}\sum_{t=1}^{T}\mbE\norm{[\H(\x_t)+\upsilon \textbf{1}]_+}
\end{align}
Thus, from \eqref{eq_zz_sum} and \eqref{eq_zz2} and using bound from \eqref{bound_on_lam2} and \eqref{divide_r2} we have:
\begin{align}\label{eq_zz3}
\sum_{t=1}^{T}\mbE[F(\x_t)-F(\x^\star_{\upsilon})]&\geq -C \sum_{t=1}^{T}\mbE\norm{[\H(\x_t)+\upsilon \textbf{1}]_+}\nonumber\\&\geq -C Q \sqrt{T}	
\end{align}
where in last inequality we used bound from \eqref{divide_r2}. Now, using Corollary \ref{coro1} and \eqref{eq_zz3} we get
\begin{align}
\sum_{t=1}^{T}\mbE[F(\x_t)&-F(\x^\star)]\geq -C Q T^{1/2} -\upsilon C T.
\label{bound_1}
\end{align}
As for \textbf{CSOA}, we have set $\eta=\frac{\hat{C}}{T^{1/2}}=\frac{C_1}{T^{1/2}}$  and $\upsilon=C_2 T^{-1/2}$ which gives:
\begin{align}
\frac{1}{T}\sum_{t=1}^{T}\mbE[F(\x_t)-F(\x^\star)]\geq& -C Q T^{-1/2}-CC_2 T^{-1/2}=-C(Q+C_2)T^{-1/2}.
\label{bound_2}
\end{align}  
where, $Q=\left[\frac{(\hat{C}^2+4\delta)(C+r)^2+D^2\hat{C}^2+5\delta C^2}{2\hat{C}r}\right]$ with $\hat{C}=C_1$ and other constants same as defined in Theorem \ref{theorem1} .

(b) As the derivation for part (a) involves working only on dual updates \eqref{iterate_dual}, which is same for both CSOA and FW-CSOA, we can obtain the lower bound on optimality gap for FW-CSOA  proceeding along the same line as in \eqref{bound_on_lam2}- \eqref{bound_2}. Thus for FW-CSOA with $\eta=\frac{\hat{C1}}{T^{3/4}}$ (which implies $\hat{C}=\hat{C1}$), we can write:
\begin{align}
\sum_{t=1}^{T}\mbE[F(\x_t)-F(\x^\star)]\geq& -C Q T^{3/4}-\upsilon CT.
\label{bound_3}
\end{align}  
As for \textbf{FW-CSOA}, we have set $\upsilon=\hat{C_2} T^{-1/4}$, we obtain:
\begin{align}
\frac{1}{T}\sum_{t=1}^{T}\mbE[F(\x_t)-F(\x^\star)]\geq& -C Q T^{-1/4}-C\hat{C_2} T^{-1/4}=-C(Q+\hat{C_2})T^{-1/4}.
\label{bound_22}
\end{align}  
where, $Q=\left[\frac{(\hat{C}^2+4\delta)(C+r)^2+D^2\hat{C}^2+5\delta C^2}{2\hat{C}r}\right]$ with $\hat{C}=\hat{C_1}$ and other constants same as defined in Theorem \ref{theorem2}.

\section{Classifier with Fairness Constraint}\label{extra_expp}

In this section we will replicate the results of \cite{zafar2015fairness}, however, instead of using batch method we will use the proposed algorithm \textbf{CSOA}. We will perform experiments on synthetic data to show that the formulation
in \eqref{app_prob_main} allows to finely control the fairness, at a small loss in
accuracy of the classifier and study the variation in $p\%$ with covariance. For a given binary sensitive attribute $s\in \{0,1\}$, for the problem defined in \eqref{app_prob_main}, we can write $p\%$ as:
\begin{align}
\label{p rule defined}
\min \left( \frac{P(\th^T\x\geq 0|s=1)}{P(\th^T\x\geq 0|s=0)},\frac{P(\th^T\x\geq 0|s=0)}{P(\th^T\x\geq 0|s=1)}\right)
\end{align}

\begin{figure}[h]
	\centering
	\setcounter{subfigure}{0}
	\begin{subfigure}{0.3\columnwidth}		\includegraphics[width=\linewidth, height = 0.6\linewidth]{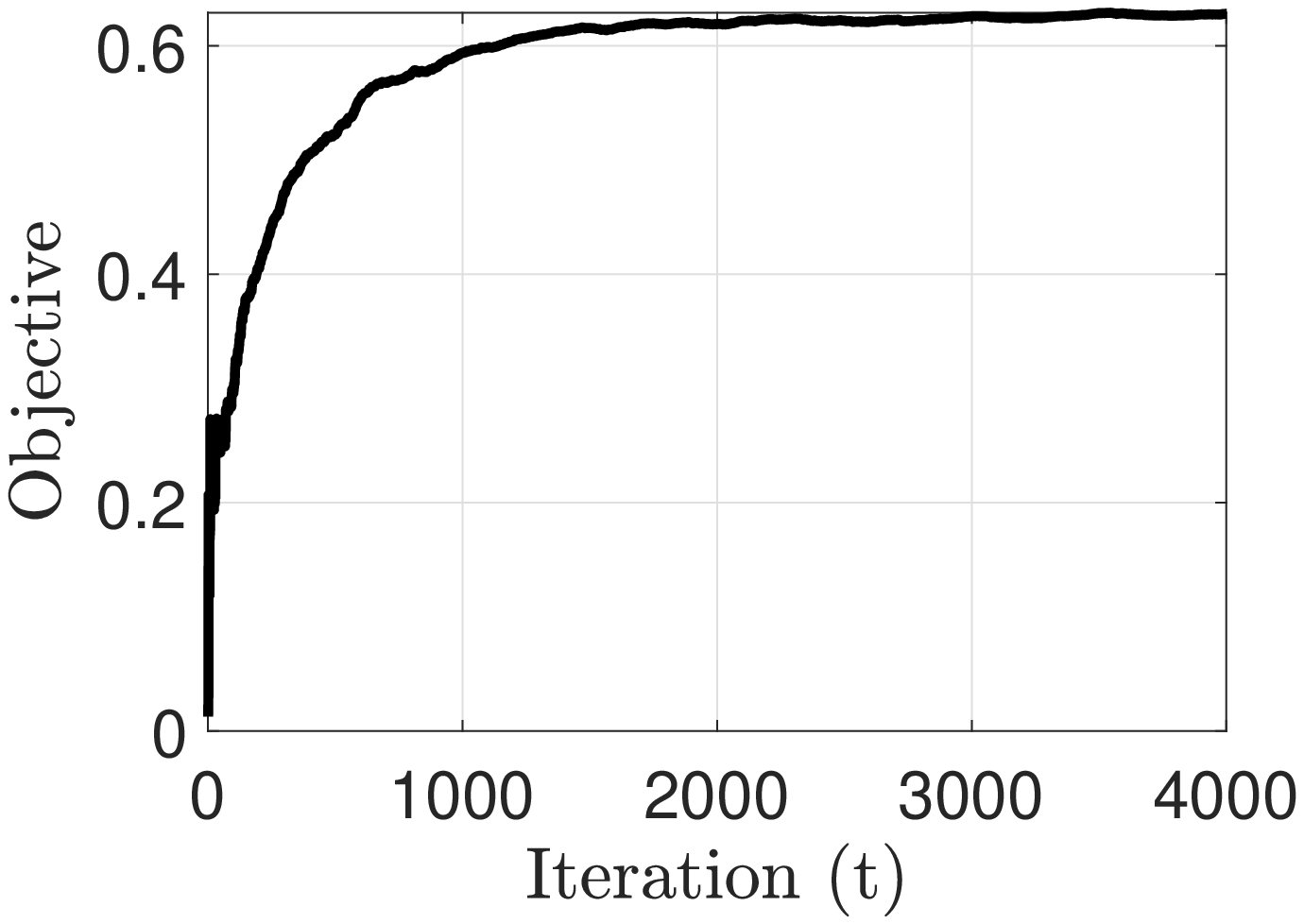}
		\caption{Normalized Error}
		\label{obj_syn}
	\end{subfigure} 	
	\begin{subfigure}{0.3\columnwidth}
		\includegraphics[width=\linewidth, height = 0.6\linewidth]{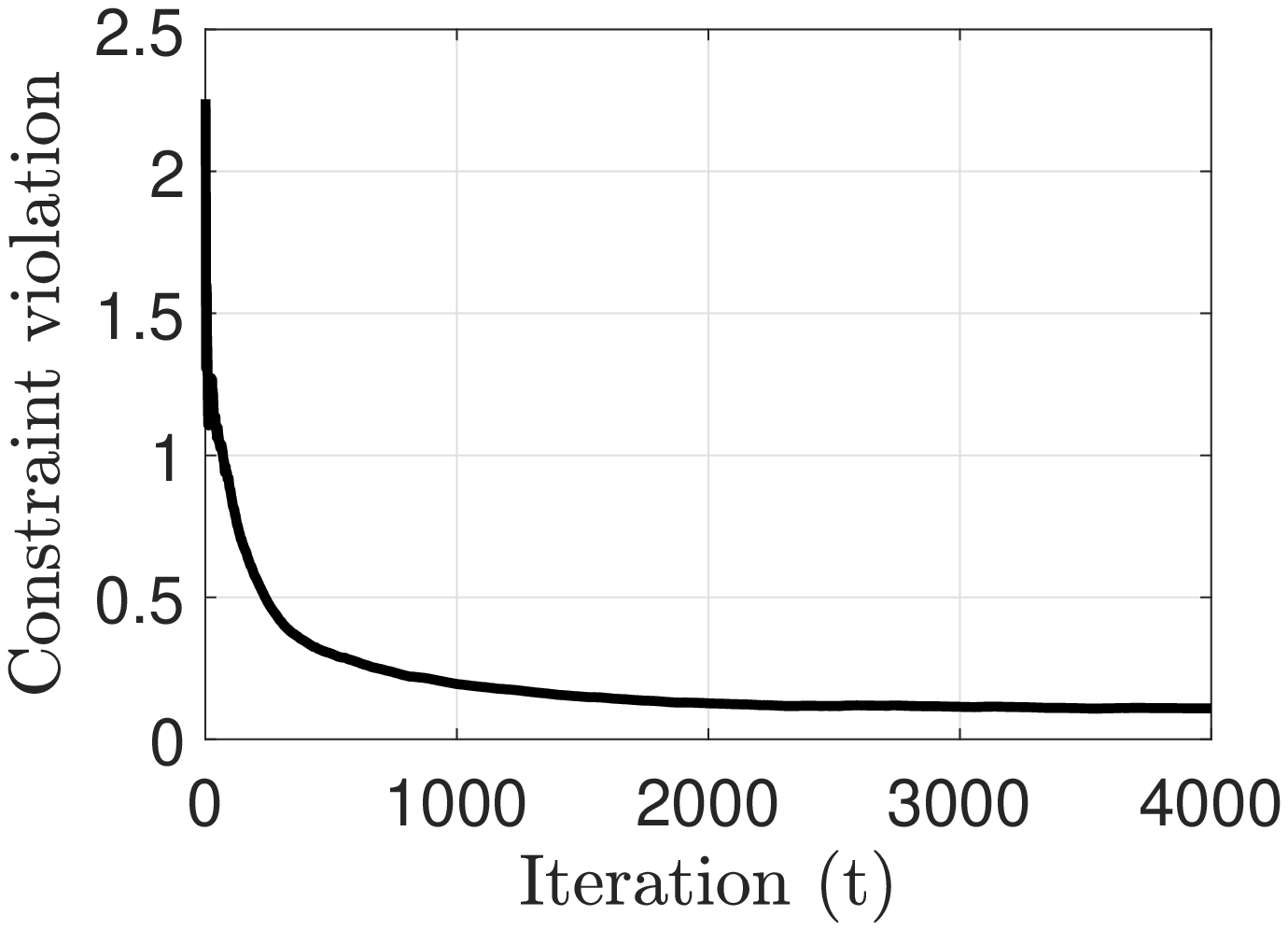}
		\caption{Constraint Violation Vs $t$ }
		\label{vio_syn}
	\end{subfigure}
	\begin{subfigure}{0.3\columnwidth}
		\includegraphics[width=\linewidth, height = 0.6\linewidth]{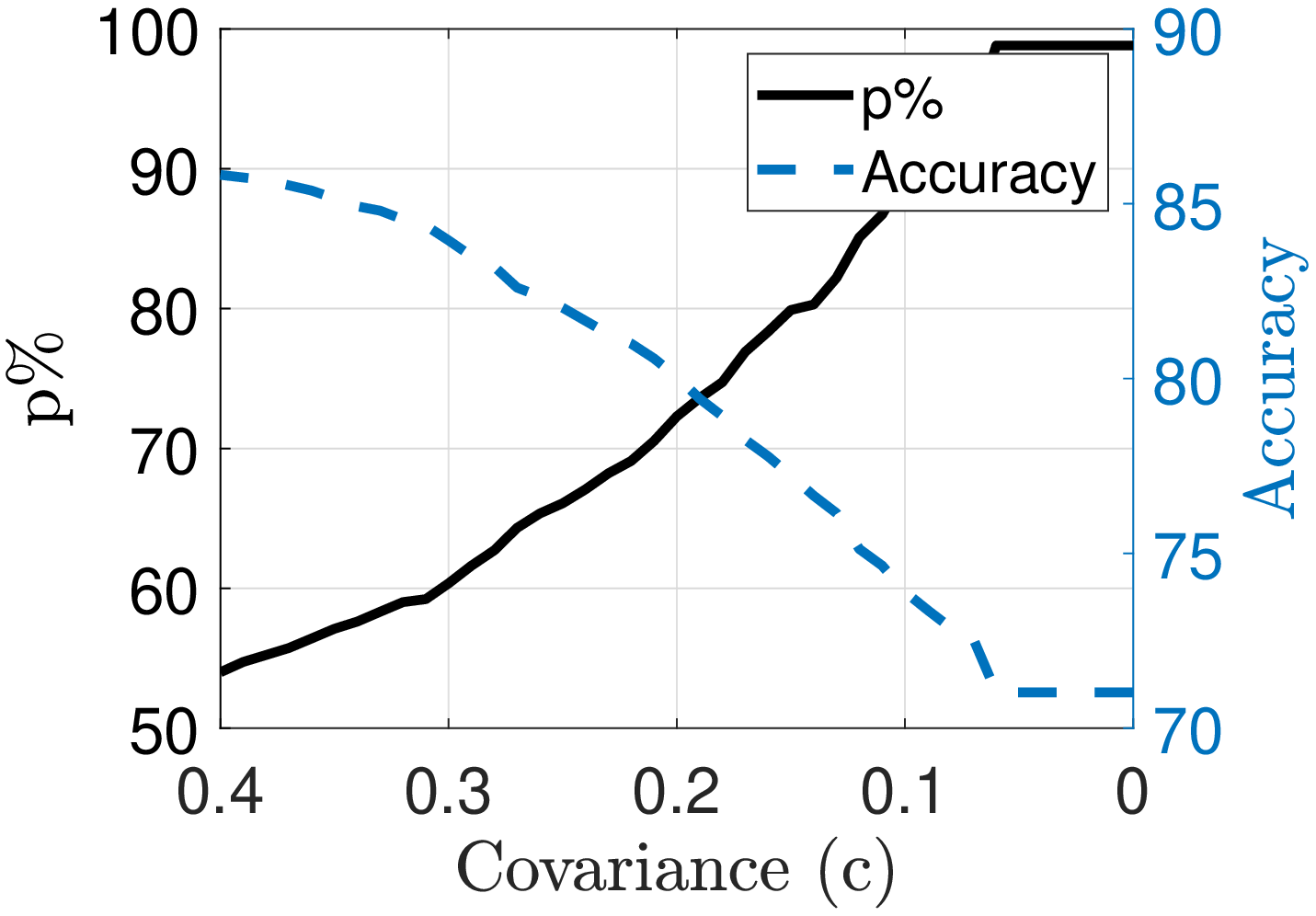}
		\caption{Variation with Covariance (C) }
		\label{P_C_SYN}
	\end{subfigure}
	\begin{subfigure}{0.25\columnwidth}
		\includegraphics[width=\linewidth, height = 0.6\linewidth]{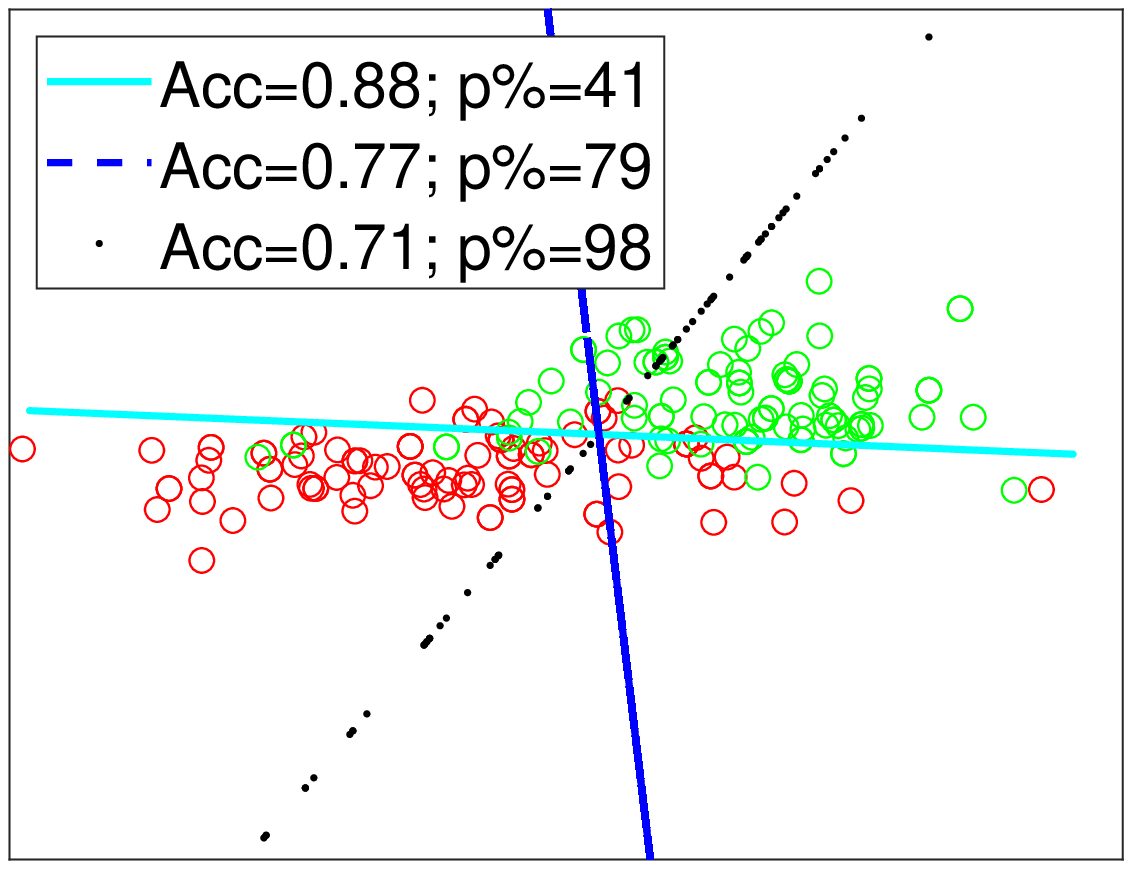}
		\caption{$\phi=\pi/4$ }
		\label{scatter_pi4}
	\end{subfigure}
	\begin{subfigure}{0.25\columnwidth}
		\includegraphics[width=\linewidth, height = 0.6\linewidth]{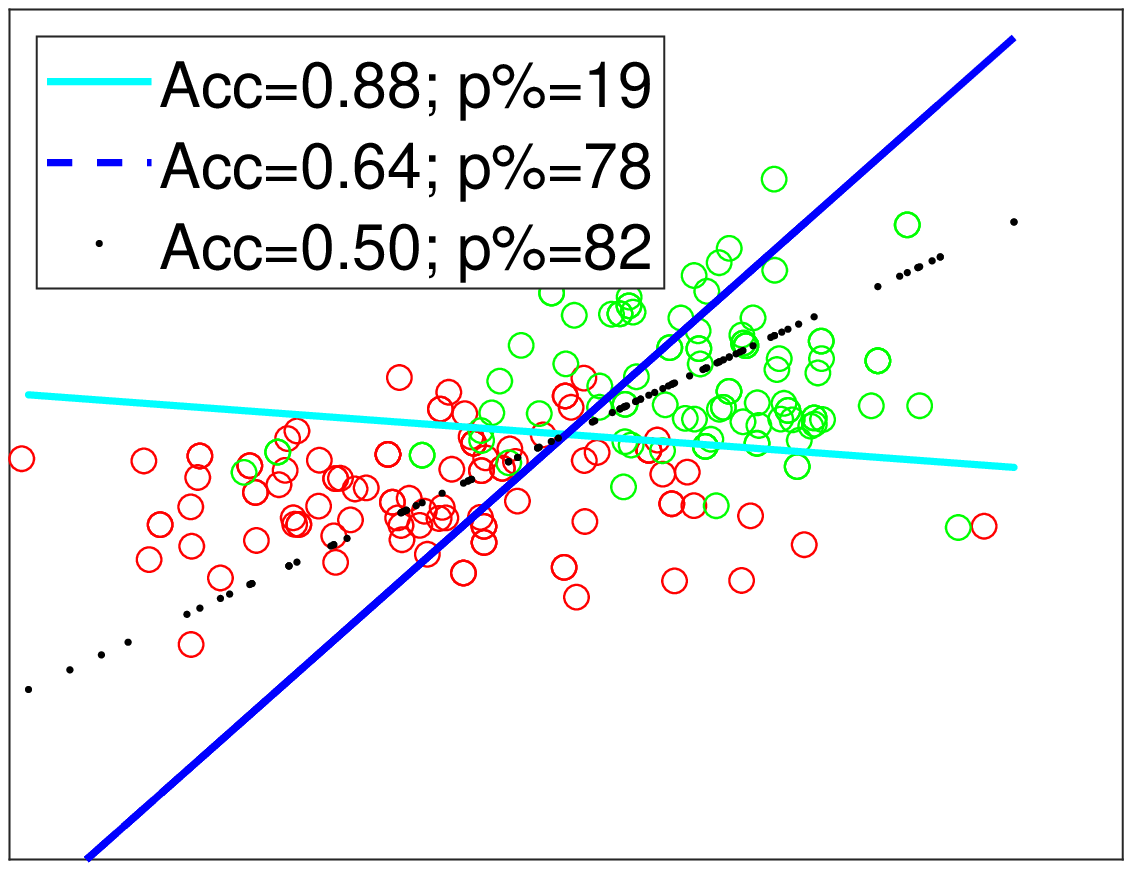}
		\caption{$\phi=\pi/8$ }
		\label{scatter_pi8}
	\end{subfigure}
	\vspace{0.25cm}
	\caption{Plots for synthetic data: (a) Objective vs Iteration, (b) Constraint violation vs Iteration. Panel (c) shows the variation in p\% and accuracy with covariance for classifiers trained under fairness constraints. For the scatter plot shown in (d) and (e), the solid light blue lines show the decision boundaries for unconstrained case for $\phi=\pi/4$ and $\phi=\pi/8$ respectively. The dashed lines show the decision boundaries trained to maximize accuracy under fairness constraints.}
	\label{fairness_syn}
\end{figure}

For the experiment, we generate synthetic data the same way as mentioned in \cite{zafar2015fairness}. Total 4000 binary labels are generated uniformly at random. Two different Gaussian distributions are generated as: $p(x|y=1)=N([2;2],[5, 1;1, 5])$ and $p(x|y=-1)=N([-2;-2],[10, 1;1, 3])$. The feature vectors are assigned to labels by drawing samples from above two distributions.  For sensitive attribute $z$ we have used Bernoulli distribution: $p(z=1)=p(\hat{\bbx}|y=1)/p(\hat{\bbx}|y=1)+p(\hat{\bbx}|y=-1)$, here $\hat{\bbx}=[\cos(\phi), -\sin(\phi);\sin(\phi), \cos(\phi)]$. The correlation between sensitive attribute and class labels (i.e. disparate impact) depends on the parameter $\phi$. The smaller the $\phi$, the higher the correlation. The problem is solved using proposed algorithm \ref{alg:CSOA}  and the results of numerical experiment are plotted in Fig. \ref{fairness_syn}.

The convergence behavior of the proposed algorithm is shown in Fig. \ref{obj_syn} and Fig. \ref{vio_syn}. It can be observed that CSOA is able to achieve almost zero constraint violation. Next, we plot the decision boundaries generated by the
classifier in Fig.\ref{scatter_pi4}-\ref{scatter_pi8}. The unconstrained
decision boundary is represented by a solid line while boundaries for different covariance thresholds $c$ is shown in dotted line. It can be observed that incorporating fairness constraint in data generation process, results into the rotation of the decision boundary (compared to unconstrained boundary). The rotation is more as the threshold value $c$ is decreased. Fig.\ref{scatter_pi4} and Fig.\ref{scatter_pi8} also shows that increase in the correlation in the data i.e. chaining $\phi$ from $\pi/4$ to $\pi/8$ causes more rotation.  Additionally, we investigate the relation between p\%, accuracy and covariance. Fig.\ref{P_C_SYN} shows that higher the covariance, the lower the p\% the classifier satisfies. Moreover, it shows the trade-off between p\% and accuracy: more p\% (i.e. smaller covariance or a larger rotation) leads to a more fair solution, but at the cost of reduced accuracy. These results shows that the problem formulation of \eqref{app_prob_main} successfully incorporates fairness constraint by exploiting the covariance and allows to control fairness constraint with small loss in accuracy.

\bibliographystyle{IEEEtran}
\bibliography{ref}

\begin{thebibliography}{10}
\providecommand{\url}[1]{#1}
\csname url@samestyle\endcsname
\providecommand{\newblock}{\relax}
\providecommand{\bibinfo}[2]{#2}
\providecommand{\BIBentrySTDinterwordspacing}{\spaceskip=0pt\relax}
\providecommand{\BIBentryALTinterwordstretchfactor}{4}
\providecommand{\BIBentryALTinterwordspacing}{\spaceskip=\fontdimen2\font plus
\BIBentryALTinterwordstretchfactor\fontdimen3\font minus
  \fontdimen4\font\relax}
\providecommand{\BIBforeignlanguage}[2]{{%
\expandafter\ifx\csname l@#1\endcsname\relax
\typeout{** WARNING: IEEEtran.bst: No hyphenation pattern has been}%
\typeout{** loaded for the language `#1'. Using the pattern for}%
\typeout{** the default language instead.}%
\else
\language=\csname l@#1\endcsname
\fi
#2}}
\providecommand{\BIBdecl}{\relax}
\BIBdecl

\bibitem{mu2017stochastic}
Y.~Mu, W.~Liu, X.~Liu, and W.~Fan, ``Stochastic gradient made stable: A
  manifold propagation approach for large-scale optimization,'' \emph{IEEE
  Trans. on Knowl. and Data Eng.}, vol.~29, no.~2, pp. 458--471, 2017.

\bibitem{li2017designing}
X.~Li, Q.~Xu, and C.~Chen, ``Designing a hierarchical decentralized system for
  distributing large-scale, cross-sector, and multipollutant control
  accountabilities,'' \emph{IEEE Systems Journal}, vol.~11, no.~4, pp.
  2774--2783, 2017.

\bibitem{derenick2007convex}
J.~C. Derenick and J.~R. Spletzer, ``Convex optimization strategies for
  coordinating large-scale robot formations,'' \emph{IEEE Trans. on Robotics},
  vol.~23, no.~6, pp. 1252--1259, 2007.

\bibitem{koppel2018parallel}
A.~Koppel, A.~Mokhtari, and A.~Ribeiro, ``Parallel stochastic successive convex
  approximation method for large-scale dictionary learning,'' in
  \emph{ICASSP}.\hskip 1em plus 0.5em minus 0.4em\relax IEEE, 2018, pp.
  2771--2775.

\bibitem{8255626}
A.~S. {Bedi} and K.~{Rajawat}, ``Network resource allocation via stochastic
  subgradient descent: Convergence rate,'' \emph{IEEE Transactions on
  Communications}, vol.~66, no.~5, pp. 2107--2121, 2018.

\bibitem{lan2012optimal}
G.~Lan, ``An optimal method for stochastic composite optimization,''
  \emph{Math. Prog.}, vol. 133, no. 1-2, pp. 365--397, 2012.

\bibitem{nemirovski2009robust}
A.~Nemirovski, A.~Juditsky, G.~Lan, and A.~Shapiro, ``Robust stochastic
  approximation approach to stochastic programming,'' \emph{SIAM Journal on
  optimization}, vol.~19, no.~4, pp. 1574--1609, 2009.

\bibitem{johnson2013accelerating}
R.~Johnson and T.~Zhang, ``Accelerating stochastic gradient descent using
  predictive variance reduction,'' in \emph{Advances in neural information
  processing systems}, 2013, pp. 315--323.

\bibitem{parikh2014proximal}
N.~Parikh, S.~Boyd \emph{et~al.}, ``Proximal algorithms,'' \emph{Found. Trends
  Optim.}, vol.~1, no.~3, pp. 127--239, 2014.

\bibitem{frank1956algorithm}
M.~Frank and P.~Wolfe, ``An algorithm for quadratic programming,'' \emph{Naval
  research logistics quarterly}, vol.~3, no. 1-2, pp. 95--110, 1956.

\bibitem{hazan2012projection}
E.~Hazan and S.~Kale, ``Projection-free online learning,'' \emph{Proceedings of
  the 29th International Conference on Machine Learning}, pp. 521--528, 2012.

\bibitem{hazan2016variance}
E.~Hazan and H.~Luo, ``Variance-reduced and projection-free stochastic
  optimization,'' in \emph{International Conference on Machine Learning}, 2016,
  pp. 1263--1271.

\bibitem{mokhtari2018stochastic}
A.~Mokhtari, H.~Hassani, and A.~Karbasi, ``Stochastic conditional gradient
  methods: From convex minimization to submodular maximization,'' \emph{Journal
  of Machine Learning Research}, vol.~21, no. 105, pp. 1--49, 2020.

\bibitem{Zee_ACC}
Z.~Akhtar and K.~Rajawat, ``Momentum based projection free stochastic
  optimization under affine constraints,'' \emph{IEEE ACC}, 2021.

\bibitem{lan2016conditional}
G.~Lan and Y.~Zhou, ``Conditional gradient sliding for convex optimization,''
  \emph{SIAM Journal on Optimization}, vol.~26, no.~2, pp. 1379--1409, 2016.

\bibitem{jaggi2013revisiting}
M.~Jaggi, ``Revisiting frank-wolfe: projection-free sparse convex
  optimization,'' in \emph{Proceedings of the 30th International Conference on
  International Conference on Machine Learning-Volume 28}, 2013, pp. I--427.

\bibitem{zafar2015fairness}
M.~B. Zafar, I.~Valera, M.~G. Rogriguez, and K.~P. Gummadi, ``Fairness
  constraints: Mechanisms for fair classification,'' in \emph{Artificial
  Intelligence and Statistics}, 2017, pp. 962--970.

\bibitem{nemirovski2006convex}
A.~Nemirovski and A.~Shapiro, ``Convex approximations of chance constrained
  programs,'' \emph{SIAM Journal on Optimization}, vol.~17, no.~4, pp.
  969--996, 2006.

\bibitem{chapelle2006semi}
O.~Chapelle, B.~Sch{\"o}lkopf, and A.~Zien, ``Semi-supervised learning, ser.
  adaptive computation and machine learning,'' 2006.

\bibitem{singh2019asynchronous}
A.~S. Bedi, A.~Koppel, and K.~Rajawat, ``Asynchronous online learning in
  multi-agent systems with proximity constraints,'' \emph{IEEE Trans. Signal
  Inf. Process. Netw.}, 2019.

\bibitem{rockafellar2000optimization}
R.~T. Rockafellar, S.~Uryasev \emph{et~al.}, ``Optimization of conditional
  value-at-risk,'' \emph{Journal of risk}, vol.~2, pp. 21--42, 2000.

\bibitem{wang2008sample}
W.~Wang and S.~Ahmed, ``Sample average approximation of expected value
  constrained stochastic programs,'' \emph{Operations Research Letters},
  vol.~36, no.~5, pp. 515--519, 2008.

\bibitem{vzliobaite2017measuring}
I.~{\v{Z}}liobait{\.e}, ``Measuring discrimination in algorithmic decision
  making,'' \emph{Data Mining and Knowledge Discovery}, vol.~31, no.~4, pp.
  1060--1089, 2017.

\bibitem{bedi2019asynchronous}
A.~S. Bedi, A.~Koppel, and K.~Rajawat, ``Asynchronous saddle point algorithm
  for stochastic optimization in heterogeneous networks,'' \emph{IEEE Trans.
  Signal Process.}, vol.~67, no.~7, pp. 1742--1757, 2019.

\bibitem{yu2017online}
H.~Yu, M.~Neely, and X.~Wei, ``Online convex optimization with stochastic
  constraints,'' in \emph{Advances in Neural Information Processing Systems},
  2017, pp. 1428--1438.

\bibitem{madavan2019subgradient}
A.~N. Madavan and S.~Bose, ``Subgradient methods for risk-sensitive
  optimization,'' \emph{arXiv preprint arXiv:1908.01086}, 2019.

\bibitem{mahdavi2012trading}
M.~Mahdavi, R.~Jin, and T.~Yang, ``Trading regret for efficiency: online convex
  optimization with long term constraints,'' \emph{J. Mach. Learn. Res.},
  vol.~13, no. Sep, pp. 2503--2528, 2012.

\bibitem{lan2020conditional}
G.~Lan and Z.~Zhou, ``Algorithms for stochastic optimization with function or
  expectation constraints,'' \emph{Computational Optimization and
  Applications}, pp. 1--38, 2020.

\bibitem{basu2019optimal}
K.~Basu and P.~Nandy, ``Optimal convergence for stochastic optimization with
  multiple expectation constraints,'' \emph{arXiv preprint arXiv:1906.03401},
  2019.

\bibitem{mahdavi2012stochastic}
M.~Mahdavi, T.~Yang, R.~Jin, S.~Zhu, and J.~Yi, ``Stochastic gradient descent
  with only one projection,'' in \emph{NeurIPS}, 2012, pp. 494--502.

\bibitem{zhang2019stochastic}
L.~Zhang, Y.~Zhang, and J.~Wu, ``Stochastic approximation proximal method of
  multipliers for convex stochastic programming,'' \emph{arXiv preprint
  arXiv:1907.12226}, 2019.

\bibitem{thomdapu2019optimal}
S.~T. Thomdapu and K.~Rajawat, ``Optimal design of queuing systems via
  compositional stochastic programming,'' \emph{IEEE Transactions on
  Communications}, vol.~67, no.~12, pp. 8460--8474, 2019.

\bibitem{thomdapu2020stochastic}
S.~T. Thomdapu, K.~Rajawat \emph{et~al.}, ``Stochastic compositional gradient
  descent under compositional constraints,'' \emph{arXiv preprint
  arXiv:2012.09400}, 2020.

\bibitem{boob2019stochastic}
D.~Boob, Q.~Deng, and G.~Lan, ``Stochastic first-order methods for convex and
  nonconvex functional constrained optimization,'' \emph{arXiv preprint
  arXiv:1908.02734}, 2019.

\bibitem{cutkosky2019momentum}
A.~Cutkosky and F.~Orabona, ``Momentum-based variance reduction in non-convex
  sgd,'' in \emph{Advances in Neural Information Processing Systems}, 2019, pp.
  15\,210--15\,219.

\bibitem{nemirovski2005efficient}
A.~Nemirovski and R.~Y. Rubinstein, ``An efficient stochastic approximation
  algorithm for stochastic saddle point problems,'' in \emph{Modeling
  Uncertainty}.\hskip 1em plus 0.5em minus 0.4em\relax Springer, 2005, pp.
  156--184.

\bibitem{shamir2013stochastic}
O.~Shamir and T.~Zhang, ``Stochastic gradient descent for non-smooth
  optimization: Convergence results and optimal averaging schemes,'' in
  \emph{International conference on machine learning}, 2013, pp. 71--79.

\bibitem{nedic2009subgradient}
A.~Nedi{\'c} and A.~Ozdaglar, ``Subgradient methods for saddle-point
  problems,'' \emph{Journal of optimization theory and applications}, vol. 142,
  no.~1, pp. 205--228, 2009.

\bibitem{asuncion20019uci}
A.~Asuncion and D.~Newman, ``Uci machine learning repository,'' 2019.

\bibitem{molitor2018matrix}
D.~Molitor and D.~Needell, ``Matrix completion for structured observations,''
  in \emph{2018 Information Theory and Applications Workshop (ITA)}.\hskip 1em
  plus 0.5em minus 0.4em\relax IEEE, 2018, pp. 1--5.

\bibitem{lan2020conditional2}
G.~Lan, E.~Romeijn, and Z.~Zhou, ``Conditional gradient methods for convex
  optimization with function constraints,'' \emph{arXiv preprint
  arXiv:2007.00153}, 2020.

\end{thebibliography}

   \end{document}